\documentclass[11pt, reqno]{amsart}
\usepackage[T1]{fontenc}
\usepackage{amsfonts}
\usepackage{amsmath}
\usepackage{amsthm}
\usepackage{amssymb}
\usepackage[textwidth=16.5cm, textheight=22cm, hmarginratio=1:1, vmarginratio=1:1]{geometry}
\usepackage{graphicx}
\usepackage{xcolor}
\usepackage{mathtools}
\usepackage[colorlinks=true, linkcolor=blue]{hyperref}
\usepackage{cleveref}
\usepackage[english]{babel}
\usepackage{lineno}
\usepackage{float}
\usepackage{booktabs,longtable,array}

\newtheorem{theorem}{Theorem}[section]

\newtheorem{lemma}[theorem]{Lemma}
\newtheorem{conjecture}[theorem]{Conjecture}
\newtheorem{proposition}[theorem]{Proposition}

\usepackage{tikz}
\usetikzlibrary{positioning}
\usetikzlibrary{decorations,arrows}
\usetikzlibrary{decorations.markings}
\numberwithin{equation}{section}

\newcommand{\E}{\mathbb E}
\newcommand{\Prb}{\mathbb P}

\newcommand{\ind}{\mathbf 1}

\newcommand{\ceil}[1]{\left\lceil#1\right\rceil}
\newcommand{\floor}[1]{\left\lfloor#1\right\rfloor}
\newcommand{\av}[1]{\left[\!\left[#1\right]\!\right]}
\DeclareMathOperator{\tr}{tr}
\DeclareMathOperator{\diag}{diag}

\IfFileExists{rustic.sty}{\usepackage{rustic}}{}
\usepackage[T1]{fontenc}

\emergencystretch=\maxdimen
\title{An improved bound on the strong chromatic index via local flag algebras}
\author{Runze Wang}
\address[]{Department of Mathematics and Computer Science, Augustana College, Rock Island, IL 61201, USA}
\email{runze.w@hotmail.com}
\thanks{}
\subjclass[2020]{05C15}

\begin{document}

\sloppy

\begin{abstract}
A strong edge coloring is an edge coloring where two edges receive different colors if they are incident or joined by a third edge. For a graph $G$, the strong chromatic index of $G$, denoted by $\chi'_s(G)$, is the minimum number of colors needed for a strong edge coloring of $G$. In this paper, using a local flag algebra method, we prove that for every graph $G$ of sufficiently large maximum degree $\Delta$, $\chi'_s(G)\le 1.6388\Delta^2$, improving the bound $1.73\Delta^2$ due to Davey, Hurley, de Joannis de Verclos, Kang, and Volec.
\end{abstract}
\keywords{}

\maketitle

\section{Introduction}

All graphs in this paper are finite and simple. A \emph{strong edge coloring} of a graph $G$ is an edge coloring where two edges receive different colors if they are incident or joined by a third edge. We say that two such edges \emph{conflict}. Equivalently, every color class is an induced matching. The \emph{strong chromatic index} of $G$, denoted by $\chi'_s(G)$, is the minimum number of colors needed for a strong edge coloring of $G$.

For graphs of bounded maximum degree, Erd\H{o}s and Ne\v{s}et\v{r}il~\cite{erdos} proposed the following conjecture.

\begin{conjecture}[Erd\H{o}s and Ne\v{s}et\v{r}il~\cite{erdos}]
\label{conj:erdos-nesetril}
Let $G$ be a graph of maximum degree $\Delta$. Then
\begin{align*}
 \chi'_s(G)\le
 \begin{cases}
   \frac{5}{4}\Delta^2 & \text{if $\Delta$ is even},\\
   \frac{5}{4}\Delta^2-\frac{1}{2}\Delta+\frac{1}{4} & \text{if $\Delta$ is odd}.
 \end{cases}
\end{align*}
\end{conjecture}

If true, the conjectured bound will be sharp, because it can be attained by a balanced blow-up of $C_5$.

By a greedy coloring, we have $\chi'_s(G)\le 2\Delta^2-2\Delta+1$. Molloy and Reed \cite{molloy} obtained the first asymptotic improvement, proving the bound $1.998\Delta^2$ for sufficiently large $\Delta$. Bruhn and Joos \cite{bruhn}, Bonamy, Perrett, and Postle \cite{bonamy}, and Hurley, de Joannis de Verclos, and Kang \cite{nibble} reduced the coefficient to $1.93$, $1.835$, and $1.772$, respectively. Very recently, Davey, Hurley, de Joannis de Verclos, Kang, and Volec \cite{strong} used \emph{local flag algebras} to obtain the bound $1.73\Delta^2$. In this paper, we reduce the coefficient to $1.6388$.

\begin{theorem}\label{thm:main}
For every graph $G$ of sufficiently large maximum degree $\Delta$,
\begin{align*}
 \chi'_s(G)\le 1.6388\Delta^2.
\end{align*}
\end{theorem}

Let $L(G)$ be the line graph of $G$ and let $X=L(G)^2$. The vertices of $X$ are the edges of $G$, and two vertices of $X$ are adjacent precisely when the corresponding edges conflict. Thus $\chi'_s(G)=\chi(X)$, where $\chi(X)$ is the chromatic number of $X$. We first delete vertices of small current degree from $X$, and seek a proper coloring of the remaining induced subgraph $Y$.

Davey et al.~\cite{strong} applied a sparse coloring theorem to an induced subgraph of $X$ with high minimum degree, using a uniform upper bound on the number of edges in each of its neighborhoods. Our local flag inequalities retain the degree of the root and a signed sum of the degrees of its neighbors. During the coloring procedure, vertices in different degree classes remain uncolored at different rates. The signed sum gives density bounds for these changing neighborhoods. Applying the local inequalities at the vertices within a neighborhood also bounds the number of triangles through each of those vertices. The inequality at the original root controls the sum of the degree terms in the triangle bounds.

The coloring procedure uses local minima of independent priorities, as in \cite{nibble}. Inclusion-exclusion gives a lower bound on the probability that a specified color is offered to no vertex in a neighborhood. We strengthen this bound by including the triangle estimates in a flag calculation on six vertices. In particular, a rooted triangle and a three-edge star with the same root involve at most six vertices. Their joint count lets us weight the triangle estimate by the cube of the degree at its root. We also use the rooted degree bound in the same calculation, with mixed moments of the rooted degree and triangle count.

The proof is computer-assisted. Appendices~\ref{app:local} and~\ref{app:ordinary} specify the flag inequalities and the rational coefficient tests, and Appendix~\ref{app:palette} gives the parameters and the recurrence for the palette bound. The flag bases, block matrices, rational witnesses, and palette schedule, together with the programs for checking the certificates and running the recurrence in exact arithmetic, are included in the supplemental materials \cite{certificates}. Readers can verify the certificates and the palette bound by running \texttt{python3 -B verify\_all.py} from the directory containing the supplemental files. The default verification uses only the standard library of Python 3.10 or later. The optional \texttt{-{}-cpp} flag replays the C++ search and compares the resulting palette bounds with the recorded file. This option requires GNU \texttt{g++} with support for C++20, OpenMP, and signed 128-bit integers.

Section~\ref{sec:prelim} gives the reductions and notation. Section~\ref{sec:density} establishes the local density inequalities, and Section~\ref{sec:offers} derives the probability estimate on which the flag calculation is based. Section~\ref{sec:profiles} incorporates the degree and triangle bounds into this estimate. Section~\ref{sec:coloring} describes the partial coloring procedure, and Section~\ref{sec:mainproof} proves Theorem~\ref{thm:main}.

\section{Preliminaries}\label{sec:prelim}

For a graph $J$, let $N_J(v)$ be the open neighborhood of $v$ and let $d_J(v)=|N_J(v)|$. For $S\subseteq V(J)$, let $J[S]$ be the subgraph induced by $S$, and let $e(J[S])=|E(J[S])|$. An ordered edge count counts the two orientations of each edge. Let $[x]_+=\max\{x,0\}$. A \emph{palette} is the collection of colors available to the coloring procedure. All logarithms are natural.

We use $o(1)$ for quantities tending to $0$ as the relevant degree scale tends to infinity. Unless otherwise stated, these errors are uniform over the graphs, edges, and neighborhoods under consideration. The number of degree classes and the number of stages of any prescribed coloring procedure are fixed before the maximum degree tends to infinity. Constants in the error estimates may depend on these fixed parameters.

\begin{lemma}\label{lem:regular}
Every graph of maximum degree at most $\Delta$ is an induced subgraph of a finite $\Delta$-regular simple graph. The conflict relation on its original edges is preserved.
\end{lemma}

\begin{proof}
Let $G_0=G$. While $G_i$ has a vertex of degree less than $\Delta$, let $G_{i+1}$ consist of two disjoint copies of $G_i$, with an edge between the two copies of each vertex of degree less than $\Delta$. Each deficient degree increases by $1$, and no degree exceeds $\Delta$. After at most $\Delta$ repetitions, the resulting graph is $\Delta$-regular. Each copy of $G_i$ is induced in $G_{i+1}$. Consequently, no edge between two original vertices is added. Whether two original edges conflict depends only on the subgraph induced by their endpoints. Since $G$ is induced in the final graph, the conflict relation is preserved.
\end{proof}

\begin{lemma}\label{lem:deletion}
Let $J$ be a graph and let $K$ be a positive integer. Repeatedly delete a vertex of current degree less than $K$, and let $J_0$ be the subgraph that remains. If $J_0$ has a proper coloring with at most $K$ colors, then so does $J$.
\end{lemma}

\begin{proof}
Reinsert the deleted vertices in reverse order. Each reinserted vertex has at most $K-1$ previously colored neighbors, so a color is available in a palette of $K$ colors.
\end{proof}

By Lemma~\ref{lem:regular}, it suffices to consider a $\Delta$-regular graph $G$. Throughout the density argument, let
\begin{align}
 X=L(G)^2,\quad D=2\Delta^2,\quad Y=X[F],\label{eq:scales}
\end{align}
where $F\subseteq E(G)$. For $f\in F$, let
\begin{align}
 d_f=\frac{d_Y(f)}{D},\quad
 c_f=\frac{2e(Y[N_Y(f)])}{D^2}.\label{eq:normalized}
\end{align}
Since $\Delta(X)\le 2\Delta(\Delta-1)$, we have $0\le c_f\le d_f^2\le1$. We refer to the degree in $Y$ as the \emph{conflict degree}. An edge $f$ under consideration is the \emph{root edge}.

\section{Local density bounds}\label{sec:density}

Let
\begin{align}
 a=0.86485,\quad b=0.6652125,\quad
 \delta=\frac{30913355784}{10^{12}},\quad k=60\delta.
 \label{eq:constants}
\end{align}
In particular, $1.8548<k<1.85481$. For a root edge $f=xy\in F$, let $B=N_G(x)\cup N_G(y)$. Let the vertices in $B$ be black and the remaining vertices be red. Let the edges in $F$ be marked. Every edge in $N_Y(f)$ has an endpoint in $B$. Let $U_1(f)$ consist of those edges with exactly one endpoint in $B$, and let $U_2(f)=N_Y(f)\setminus U_1(f)$.

\begin{lemma}\label{lem:incidence}
For every $f\in F$,
\begin{align}
 \frac{|U_2(f)|}{D}\le 1-d_f.\label{eq:incidence}
\end{align}
\end{lemma}

\begin{proof}
There are at most $|B|\Delta\le2\Delta^2=D$ incidences at vertices in $B$. Each edge in $U_1(f)$ contributes one incidence, and each edge in $U_2(f)$ contributes two. Thus $|U_1(f)|+2|U_2(f)|\le D$. Since $|U_1(f)|+|U_2(f)|=d_fD$, the conclusion follows.
\end{proof}

Let $\theta=0.8175$.
\begin{lemma}\label{lem:density}
For every $F\subseteq E(G)$ and every $f\in F$,
\begin{align}
 c_f\le b+\frac{k}{D}\sum_{e\in U_1(f)}(a-d_e)+o(1).
 \label{eq:density-base}
\end{align}
Assume in addition that $d_e\ge\theta-\frac{1}{D}$ for every $e\in F$. For each pair $(\gamma,b_\gamma)$ in Table~\ref{tab:root},
\begin{align}
 c_f\le b_\gamma+\frac{k}{D}\sum_{e\in U_1(f)}(a-d_e)
           -\gamma(1-d_f)+o(1).\label{eq:density-root}
\end{align}
For each triple $(\kappa,\gamma,b')$ in Table~\ref{tab:extra-root}, we also have
\begin{align}
 c_f\le b'-\gamma(1-d_f)+\frac{\kappa}{D}\sum_{e\in U_1(f)}(a-d_e)+o(1).
 \label{eq:density-extra}
\end{align}
The errors are uniform as $\Delta\to\infty$.
\end{lemma}

\begin{table}[ht]
\centering
\caption{Parameters in the inequalities involving the root degree.}\label{tab:root}
\begin{tabular}{cc}
\toprule
$\gamma$&$b_\gamma$\\
\midrule
0.38&0.66577886\\
0.50&0.66713714\\
0.60&0.66776249\\
0.64&0.66967403\\
0.665&0.67078806\\
\bottomrule
\end{tabular}
\end{table}

\begin{table}[htbp]
\centering
\caption{Additional parameters in the inequalities involving the root degree.}\label{tab:extra-root}
\begin{tabular}{@{}ccc@{}}
\toprule
$\kappa$ & $\gamma$ & $b'$\\
\midrule
$0$ & $0.73657442$ & $0.73657943$\\
$0$ & $0$ & $0.72861074$\\
$\frac{4}{5}-\frac{1}{5\cdot 10^{10}}$ & $0$ & $0.69104375$\\
$1.2$ & $0$ & $0.67550863$\\
\bottomrule
\end{tabular}
\end{table}

\subsection{The signed degree expression}

We use the local flag normalization of \cite{local,strong}. A colored graph is local if each of its connected components contains a black vertex. For a local unlabeled colored graph $H$ on $r$ vertices, let its density be the number of induced copies of $H$ divided by $\binom{\Delta}{r}$. Marks on edges form part of the induced colored graph. Appendix~\ref{app:local} gives the product, averaging, and coefficient conventions used below.

Let $O$ be the size-$5$ expression that counts two disjoint conflicting marked edges, each with a black endpoint, with an additional adjacent vertex. Let $S$ be the signed size-$5$ expression for the conflict degree of a marked red-black edge. The precise coefficients are given in \eqref{eq:flag-objective}--\eqref{eq:flag-signed}. Their evaluations on the colored graph associated with $f$ satisfy
\begin{align}
 O&=16c_f+o(1),\label{eq:bridge-o}\\
 S&=\frac{240}{\Delta^2}\sum_{e\in U_1(f)}
       \left(\frac{d_Y(e)}{\Delta^2}-1.7297\right)+o(1).
 \label{eq:bridge-s}
\end{align}
We verify \eqref{eq:bridge-s} by counting two types of extensions. Let $p$ count a marked red-black edge, a disjoint conflicting marked edge, and one further neighbor of the red endpoint of the first edge. Let $q$ count a marked red-black edge and an ordered triple of distinct further neighbors of its red endpoint. The size-$5$ coefficient is $2p-3.4594q$. For a fixed first edge $e$, the two counts are respectively $\Delta d_Y(e)+O(\Delta^2)$ and $\Delta^3+O(\Delta^2)$. Division by $\binom{\Delta}{5}$ gives \eqref{eq:bridge-s}. For \eqref{eq:bridge-o}, each unordered pair of nonincident conflicting edges in $N_Y(f)$ contributes $4\Delta+O(1)$ to the count in \eqref{eq:o-coeff}. There are also $O(\Delta^2)$ pairs containing the root edge $f$, each contributing $O(\Delta)$. Their contribution after division by $60\binom{\Delta}{5}$ is $O(\Delta^{-2})$. There are $O(\Delta^3)$ incident pairs in $Y[N_Y(f)]$, so they contribute $O(\Delta^{-1})$ to $c_f$. Dividing the objective count by $60\binom{\Delta}{5}$ now gives \eqref{eq:bridge-o}.

The calculation in Proposition~\ref{prop:local-finite} gives
\begin{align}
 O+\delta S\le10.6434+o(1).\label{eq:signed-local}
\end{align}
The nonnegative expressions used in this calculation come from squares of flags, regularity of $G$, and the bound $|B|\le2\Delta$. They apply for every choice of $F$. Since $\frac{10.6434}{16}=b$ and $k=60\delta$, substitution of \eqref{eq:bridge-o} and \eqref{eq:bridge-s} proves \eqref{eq:density-base}.

\subsection{The degree of the root edge}

Let $F^+$ consist of $F$ together with every edge whose endpoints are both red. Mark the edges of $F^+$. The marked edges with an endpoint in $B$ remain the same, so the neighborhood of $f$, $c_f$, and $U_1(f)$ remain the same. Let $Y^+=X[F^+]$. For every marked edge $e$ with a black endpoint, we have $d_{Y^+}(e)\ge d_Y(e)\ge\theta D-1$.

Let $m$ be the number of marked edges with an endpoint in $B$, divided by $D$. The root edge itself is one of these edges, and the others form $N_Y(f)$. Thus
\begin{align}
 m=d_f+\frac{1}{D}.\label{eq:root-mass}
\end{align}
Equation~\eqref{eq:flag-mass} gives the corresponding local flag expression up to $o(1)$.

For each ordered marked edge $xy$ with a black endpoint, let $C=N_G(x)\cup N_G(y)$. At least $\theta D-1$ marked edges have an endpoint in $C$. In particular,
\begin{align}
 |N_G(x)\cap N_G(y)|\le 2(1-\theta)\Delta+O(1).
 \label{eq:codegree}
\end{align}
The count $\Delta|C|-e(G[C])$ of all edges with an endpoint in $C$ gives a stronger inequality when we also account for the unmarked edges incident with $C$. Appendix~\ref{app:local} uses these inequalities after conditioning on local flags. It also uses the lower bound on $m$ and bounds on the number of black vertices and the number of marked incidences at those vertices.

By Proposition~\ref{prop:local-finite}, for each row of Table~\ref{tab:root},
\begin{align}
 O+\delta S-16\gamma m\le16(b_\gamma-\gamma)+o(1),
 \label{eq:root-local}
\end{align}
where $S$ is evaluated with the marks $F^+$. Equations \eqref{eq:bridge-o}--\eqref{eq:bridge-s} and \eqref{eq:root-mass} give the required bound with the conflict degrees in $Y^+$. Replacing those degrees by the degrees in $Y$ increases the right-hand side and proves \eqref{eq:density-root}. For a triple $(\kappa,\gamma,b')$ in Table~\ref{tab:extra-root}, Proposition~\ref{prop:local-finite} similarly gives
\begin{align*}
 O+\frac{\kappa}{60}S-16\gamma m\le 16(b'-\gamma)+o(1).
\end{align*}
The same substitution proves \eqref{eq:density-extra}, since $\kappa\ge0$.

\subsection{A bound for local flag densities}

The rational verification of \eqref{eq:root-local} permits a small coefficient error. The following bound controls its total contribution.

\begin{lemma}\label{lem:tree}
For a local colored graph $H$ on five vertices, let $h(H)$ be the number of spanning trees of the graph obtained by adding a vertex adjacent to all black vertices of $H$. If $\varphi(H)$ is its local density, then
\begin{align}
 \sum_H h(H)\varphi(H)\le4802+o(1).\label{eq:tree-bound}
\end{align}
The sum ranges over the local size-$5$ colored graphs that occur in the marked graph under consideration.
\end{lemma}

\begin{proof}
Label the five vertices by $1,\ldots,5$ and the added vertex by $0$. For a spanning tree $T$, each child of $0$ has at most $2\Delta$ choices, and every other vertex has at most $\Delta$ choices once its parent has been chosen. Thus $T$ has at most $2^{d_T(0)}\Delta^5$ such embeddings. The weighted form of Cayley's formula, or the Pr\"ufer code, gives
\begin{align*}
 \sum_T2^{d_T(0)}=2(2+5)^4=4802.
\end{align*}
Each unlabeled five-vertex subset has $5!$ labelings. Division by $5!\binom{\Delta}{5}=\Delta^5+O(\Delta^4)$ proves \eqref{eq:tree-bound}. Every local graph has $h(H)\ge1$.
\end{proof}

\section{Color-offer probabilities}\label{sec:offers}

\subsection{The random independent set}

Let $J$ be a graph with $\Delta(J)\le M$, where $M$ is a positive integer. For each $v\in V(J)$, let $\xi_v$ be uniform on $[0,1]$. For each $v$, also let $\zeta_{v,1},\ldots,\zeta_{v,M-d_J(v)}$ be uniform on $[0,1]$. All these variables are independent. Let
\begin{align}
 \mathcal B(v)=\{\xi_u:u\in N_J(v)\}
       \cup\{\zeta_{v,i}:1\le i\le M-d_J(v)\}.
 \label{eq:comparisons}
\end{align}
Thus $\mathcal B(v)$ consists of exactly $M$ comparison variables. The variables $\zeta_{v,i}$ are used only in comparisons at $v$. Let $I$ consist of the vertices $v$ for which $\xi_v$ is less than every variable in $\mathcal B(v)$. We say that these vertices receive an offer. Ties have probability zero.

\begin{lemma}\label{lem:priority}
The set $I$ is independent, and $\Prb(v\in I)=\frac{1}{M+1}$ for every $v\in V(J)$. For distinct vertices $u,v$,
\begin{align}
 \Prb(u,v\in I)\le\frac{2}{(M+1)(M+2)}.\label{eq:pair-universal}
\end{align}
\end{lemma}

\begin{proof}
If $uv\in E(J)$, membership of $u$ in $I$ requires $\xi_u<\xi_v$, and membership of $v$ requires the reverse inequality. Hence $I$ is independent. The $M+1$ variables at a fixed vertex are independent and identically distributed, which gives the one-vertex probability.

For nonadjacent $u,v$, let $U=|\mathcal B(u)\cup\mathcal B(v)|$. The integration in \eqref{eq:joint-integral} with the original comparison sets gives $\Prb(u,v\in I)=\frac{2}{(M+1)(U+2)}$. Since $U\ge M$, this proves \eqref{eq:pair-universal}.
\end{proof}

Let $S\subseteq V(J)$ with $|S|\le M$, and let
\begin{align}
 t=\frac{|S|}{M},\quad c=\frac{2e(J[S])}{M^2}.\label{eq:tc}
\end{align}
Then $0\le c\le t^2\le1$. Let
\begin{align}
 g(s)=\frac{s}{2}-\frac{s^{\frac{3}{2}}}{6}.\label{eq:g}
\end{align}

\subsection{Independent comparisons outside the test set}

We use a reduction of Hurley, de Joannis de Verclos, and Kang \cite[Claim~2 in the proof of Theorem~2.1]{nibble}. For each priority of a vertex outside $S$ that appears in more than one comparison set at vertices of $S$, replace that priority by an independent copy in each such set. All comparison variables other than priorities of vertices of $S$ are then private to a single candidate. For $v\in S$, let $\mathcal B_{\rm split}(v)$ be its comparison set after the replacement. Let $I_{\rm split}$ consist of the vertices of $S$ that receive offers in this model.

\begin{lemma}\label{lem:split}
We have $\Prb(I\cap S\ne\varnothing)\le \Prb(I_{\rm split}\ne\varnothing)$.
\end{lemma}

\begin{proof}
Consider one shared variable and condition on all other variables. Let $A$ be the vertices of $S$ that would receive an offer in the absence of this comparison, and that use it. If another vertex already receives an offer regardless of the shared variable, the desired comparison is immediate. Otherwise, let $v\in A$ have the smallest priority, provided $A$ is nonempty. Give $v$ the original shared variable and give independent copies to all other vertices using it. An offer to any vertex of $A$ before the replacement implies that the original variable exceeds $\xi_v$, so $v$ receives an offer after the replacement. Conditional on the other variables, the copies have the required independent uniform distribution. Repeat this operation for all shared variables outside $S$.
\end{proof}

Two vertices of $S$ share comparison variables in the split model only through their common neighbors in $J[S]$. Let $Q_S$ be the complement of $J[S]$. For $u,v\in S$, let
\begin{align}
 b_{uv}=\frac{|N_{J[S]}(u)\cap N_{J[S]}(v)|}{M}.
 \label{eq:buv}
\end{align}
A clique in $Q_S$ is an independent set in $J[S]$.

Let $A$ be an independent $j$-set in $J[S]$, where $j$ is fixed. For an ordering $\pi$ of its vertices, let
\begin{align*}
 U_i(\pi)=\left|\bigcup_{h=1}^i\mathcal B_{\rm split}(\pi(h))\right|.
\end{align*}
The comparison variables are identified by the underlying random variable, so a shared variable is counted once in this union. An exact integration gives
\begin{align}
 M^j\Prb(A\subseteq I_{\rm split})
   =\sum_{\pi}\prod_{i=1}^j
      \left(\frac{U_i(\pi)}{M}+\frac{i}{M}\right)^{-1}.
 \label{eq:joint-integral}
\end{align}
To verify this identity, let $T_v=-M\log(1-\xi_v)$ and apply the same transformation to every comparison variable. On the region where the $T_v$ are in decreasing order $\pi$, integrate over the successive nonnegative differences. The exponential decay rate in the $i$th difference is $\frac{U_i(\pi)+i}{M}$. The identity holds for the original model with the corresponding comparison sets. Since the true union sizes satisfy $U_i(\pi)\ge M$, omitting $\frac{i}{M}$ in each factor gives an approximation with error $O(M^{-1})$ in the scaled probability, uniformly for fixed $j$. We call the resulting expression the leading weight of $A$.

\subsection{Pairs and triples}

We use the pair and triple argument from the proof of \cite[Theorem~2.1]{nibble}. We first replace the scaled probabilities in \eqref{eq:joint-integral} by the leading weights. In the inclusion-exclusion sum for $j=2,3$, the total error is $O(M^{-1})$, since there are at most $\binom{M}{j}$ terms and each probability has error $O(M^{-j-1})$. For a nonadjacent pair in $J[S]$, its leading weight in \eqref{eq:joint-integral} is
\begin{align*}
 \frac{2}{2-b_{uv}}.
\end{align*}
For an independent triple, let $b_1,b_2,b_3$ be its three pair overlaps. Each is at most $\frac{M-2}{M}$, since the two vertices of a pair lie outside their common neighborhood. The leading weight of the triple is at most
\begin{align}
 \frac{2\sum_{i=1}^3\frac{1}{2-b_i}}{3-\sum_{i=1}^3 b_i}
 \le 1+\sum_{i=1}^3
       \frac{b_i(3-b_i)}{3(2-b_i)(1-b_i)}.
 \label{eq:triple-comparison}
\end{align}
For the second inequality, convexity gives $\frac{1}{1-\overline b}\le\frac{1}{3}\sum_i\frac{1}{1-b_i}$, where $\overline b=\frac{1}{3}\sum_i b_i$. Chebyshev's sum inequality for the increasing functions $\frac{1}{2-b}$ and $\frac{1}{1-b}$ then gives
\begin{align*}
 \frac{2}{3}\left(\sum_i\frac{1}{2-b_i}\right)
      \frac{1}{1-\overline b}
 &\le\frac{2}{3}\sum_i\frac{1}{(2-b_i)(1-b_i)}\\
 &=1+\sum_i\frac{b_i(3-b_i)}{3(2-b_i)(1-b_i)}.
\end{align*}
The first inequality in \eqref{eq:triple-comparison} follows by bounding the union of three comparison sets from below by $M(3-b_1-b_2-b_3)$.

A pair with overlap $b$ extends to at most $(1-b)M$ independent triples in $S$, since an extension cannot be a common neighbor of the pair. Charge the term for this pair in \eqref{eq:triple-comparison} to each independent triple containing it. After multiplying the pair and triple terms in the inclusion-exclusion sum by $M^{-2}$ and $M^{-3}$, respectively, the pair weight above $1$ less the total charge is at least $M^{-2}$ times
\begin{align*}
 \frac{b}{2-b}-\frac{b(3-b)}{3(2-b)}
   =\frac{b^2}{3(2-b)}\ge0.
\end{align*}
Let $s=t^2-c$. Since $Q_S$ is the complement of $J[S]$, we have $\frac{2e(Q_S)}{M^2}=s-\frac{t}{M}\le s$. If $A_Q$ is the adjacency matrix of $Q_S$, then the number of triangles in $Q_S$ is $\frac{1}{6}\tr(A_Q^3)$. Since $\tr(A_Q^3)\le(\tr(A_Q^2))^{\frac{3}{2}}$ and $\tr(A_Q^2)=2e(Q_S)$, this number is at most $\frac{1}{6}M^3s^{\frac{3}{2}}$. Three terms of inclusion-exclusion and the preceding comparison bound the probability of an offer in the split model. Thus
\begin{align}
 \Prb(I\cap S\ne\varnothing)\le\Prb(I_{\rm split}\ne\varnothing)
   \le t-\frac{s}{2}+\frac{s^{\frac{3}{2}}}{6}+o(1).
 \label{eq:three-term}
\end{align}
The proof of the second inequality applies to any set of at most $M$ candidates whose outside comparison variables are private.

\subsection{Enlarging the test set}

Add $M-|S|$ auxiliary candidates to the split model, and let $S^+$ be the enlarged candidate set. Let $H$ be the graph on $S^+$ obtained from $J[S]$ by adding these candidates as isolated vertices. Each auxiliary candidate has its own priority and $M$ independent private comparison variables. For $v\in S^+$, let $\mathcal B_{\rm pad}(v)$ be its comparison set, which is $\mathcal B_{\rm split}(v)$ if $v\in S$. Let $I_{\rm pad}$ be the set of candidates receiving offers. Offers to the auxiliary candidates are independent of offers to the original vertices, so
\begin{align}
 \Prb(I\cap S=\varnothing)\ge\Prb(I_{\rm split}=\varnothing)
 =\frac{\Prb(I_{\rm pad}=\varnothing)}
         {\left(1-\frac{1}{M+1}\right)^{M-|S|}}.
 \label{eq:padding}
\end{align}
The denominator is $e^{t-1}+o(1)$ uniformly in $t$. Let $Q=\overline H$ and let $Z=|I_{\rm pad}|$. Then $\frac{2e(Q)}{M^2}=1-c-\frac{1}{M}$. Every candidate receives an offer with probability $\frac{1}{M+1}$, so $\E Z=1+O(M^{-1})$. The identity \eqref{eq:joint-integral} applies with $I_{\rm pad}$ and $\mathcal B_{\rm pad}$ in this model.

\begin{lemma}\label{lem:four-five}
The numbers of offers satisfy
\begin{align}
 \E\binom Z5\le\frac{1}{5}\E\binom Z4.\label{eq:five-four}
\end{align}
Consequently,
\begin{align}
 \Prb(Z=0)\ge\E\binom Z2-\E\binom Z3
                       +\frac{4}{5}\E\binom Z4+o(1). \label{eq:bonferroni}
\end{align}
\end{lemma}

\begin{proof}
Let $A$ be an independent $4$-set in $H$, and condition on $A\subseteq I_{\rm pad}$. Let $m_A$ be the number of candidate priorities in the union of the comparison sets at the vertices of $A$. At most $M-m_A$ candidates can extend $A$ to an independent $5$-set. Such a candidate $v$ has no neighbor in $A$, so its priority is independent of $A\subseteq I_{\rm pad}$ and none of its comparisons uses a priority of a vertex of $A$. Every outside comparison is private, and at most $m_A$ of the comparisons at $v$ can be shared with those at $A$. At least $M-m_A$ comparisons at $v$ and the priority of $v$ are therefore independent of the conditioning event. Thus
\begin{align*}
 \Prb(v\in I_{\rm pad}\mid A\subseteq I_{\rm pad})
 \le\frac{1}{M-m_A+1}.
\end{align*}
Sum over the possible extensions and then over $A$. Each surviving $5$-set is counted five times, which proves \eqref{eq:five-four}. Bonferroni's inequality gives
\begin{align*}
 \Prb(Z=0)\ge1-\E Z+\E\binom Z2-\E\binom Z3
                     +\E\binom Z4-\E\binom Z5.
\end{align*}
Now use $\E Z=1+O(M^{-1})$ and \eqref{eq:five-four}.
\end{proof}

\subsection{A polynomial on six vertices}

All expectations in the following calculation refer to independent uniform samples from the $M$ candidates. Coincident samples affect the expressions by $O(M^{-1})$. Let $\mathcal K_j$ be the event that the first $j$ sampled vertices form a clique in $Q$. For a pair of sampled vertices, let $b_{12}$ be the normalized number of common neighbors in the graph complementary to $Q$. For three sampled vertices, let $b_{123}$ be their normalized common-neighbor count in that graph. On $\mathcal K_3$, let
\begin{align}
 x=\frac{b_{12}}{2},\quad
 y=\frac{b_{12}+b_{13}+b_{23}-b_{123}}{3}.
 \label{eq:xy}
\end{align}
The union of three comparison sets has normalized size at least $1$, so $0\le x\le\frac{1}{2}$ and $0\le y\le\frac{2}{3}$. Let
\begin{align}
 \Psi(x,y)=\sum_{\substack{i,j\ge0\\i+j<3}}x^iy^j
                 +6x^3+3x^2y+3xy^2+3y^3.\label{eq:Psi}
\end{align}
For these $x,y$,
\begin{align}
 \frac{1}{(1-x)(1-y)}\le\Psi(x,y).\label{eq:geom-bound}
\end{align}
Indeed, the remainder after the terms of total degree less than $3$ is
\begin{align*}
 \frac{x^3}{(1-x)(1-y)}
        +\frac{x^2y+xy^2+y^3}{1-y},
\end{align*}
which is at most $6x^3+3(x^2y+xy^2+y^3)$.

On $\mathcal K_4$, let $\overline b$ be the mean of the six pair overlaps. For an ordering of the four vertices, subtract the intersections with the first comparison set to obtain
\begin{align*}
 \left|\bigcup_{h=1}^i\mathcal B_{\rm pad}(\pi(h))\right|
 \le iM-\sum_{h=2}^i
      |\mathcal B_{\rm pad}(\pi(1))\cap\mathcal B_{\rm pad}(\pi(h))|.
\end{align*}
Averaging over the orderings bounds the mean normalized union size by $i-(i-1)\overline b$. The product of reciprocals is convex on the positive orthant. Jensen's inequality in \eqref{eq:joint-integral} therefore bounds the leading weight of the $4$-set from below by
\begin{align}
 \prod_{a\in\{\frac{1}{2},\frac{2}{3},\frac{3}{4}\}}
     \frac{1}{1-a\overline b}.
 \label{eq:four-weight}
\end{align}
Let
\begin{align}
 h_1=\frac{1}{2}+\frac{2}{3}+\frac{3}{4}=\frac{23}{12},\quad
 h_2=\sum_{a\in A}a^2+\sum_{\substack{a,a'\in A\\a<a'}}aa',
       \quad A=\left\{\frac{1}{2},\frac{2}{3},\frac{3}{4}\right\}.
 \label{eq:h-coeff}
\end{align}
The expression in \eqref{eq:four-weight} is at least $1+h_1\overline b+h_2\overline b^2$.

Define $\mathcal F(Q)$ by
\begin{align}
 \mathcal F(Q)=
 \frac{1}{2}\E\left[\ind_{\mathcal K_2}\sum_{i=0}^4\left(\frac{b_{12}}{2}\right)^i\right]
 -\frac{1}{6}\E\left[\ind_{\mathcal K_3}\Psi(x,y)\right]
 +\frac{1}{30}\E\left[\ind_{\mathcal K_4}(1+h_1\overline b+h_2\overline b^2)\right].
 \label{eq:F}
\end{align}
The geometric expansion of the pair weight gives the first term. Equation~\eqref{eq:geom-bound}, together with the symmetry of the ordering in \eqref{eq:joint-integral}, bounds the triple term from above. Equation~\eqref{eq:four-weight} bounds the $4$-set term from below. Thus \eqref{eq:bonferroni} gives
\begin{align}
 \Prb(Z=0)\ge\mathcal F(Q)-o(1).\label{eq:no-offer-F}
\end{align}

Up to an $O(M^{-1})$ error, $\mathcal F(Q)$ is a linear combination of induced densities of graphs on six vertices. Each factor $b_{12}$ uses one additional sampled vertex. For a fixed triangle of $Q$, the expression $3y$ is the expectation of $\max\{0,2-j\}$, where $j$ is the number of neighbors in that triangle of one new sampled vertex. For a fixed $4$-clique, the expression $6\overline b$ is the expectation of $\binom{4-j}{2}$ for one new sampled vertex. Use a separate independent draw for each factor in a product of these expectations. The largest number of samples needed in \eqref{eq:F} is six.

For $s\in[0.25,1]$, let $\beta(s)$ be the function specified in \eqref{eq:bonus} using the bonuses in Table~\ref{tab:ordinary}. Let $T_1$ and $T_2$ be the two nondecreasing functions specified in Table~\ref{tab:tail-gains}. Let
\begin{align}
 L_0(s)=\max\bigl\{
 0.018+0.385\min\{s,0.446\}+\beta(s),\ T_1(s),\ T_2(s),\ g(s)
 \bigr\}.
 \label{eq:L0}
\end{align}
Proposition~\ref{prop:ordinary-finite} gives $\Prb(Z=0)\ge L_0(1-c)-o(1)$ whenever $1-c\ge0.25$. Each of the four functions in \eqref{eq:L0} is nondecreasing. At $s=0.43$, the first term in \eqref{eq:L0} is $0.19096$, which exceeds both $T_1(0.43)=0.18960$ and $T_2(0.43)=0.19090$, so $L_0$ is continuous there. Let $c_+\in[0,0.75]$, and assume $c\le c_++o(1)$. By monotonicity and continuity of $L_0$ and the endpoint estimate in Appendix~\ref{app:ordinary}, we have $\Prb(Z=0)\ge L_0(1-c_+)-o(1)$. If $t_+$ is an upper bound on $t$, then \eqref{eq:padding} gives
\begin{align}
 \Prb(I\cap S\ne\varnothing)
 \le 1-L_0(1-c_+)e^{1-t_+}+o(1).
 \label{eq:baseline-offer}
\end{align}
The same argument applies to any other nonnegative lower bound for the probability of no offer in the enlarged set. We obtain such bounds in Section~\ref{sec:profiles}.

\section{Weighted neighborhood bounds}\label{sec:profiles}

Let $\tau=\theta=0.8175$, let $\tau_1=0.8194$, and let $K=\floor{\tau_1D}$. Let $Y=X[F]$ be the induced subgraph remaining after repeatedly deleting vertices of current degree less than $K$. Since $\tau_1\ge\tau$, every edge of $F$ satisfies $d_e\ge\tau-\frac{1}{D}$. We may assume that $F\ne\varnothing$. For $e\in F$, let $\widehat d_e=\max\{\tau,d_e\}$ and let $x_e=\widehat d_e-\tau$. Thus $0\le \widehat d_e-d_e\le\frac{1}{D}$. The inequalities of Section~\ref{sec:density} hold with $\widehat d_e$ in place of $d_e$, with the same uniform $o(1)$ error. Let $\mathcal P$ consist of $(0,b)$ and the five pairs $(\gamma,b_\gamma)$ in Table~\ref{tab:root}. Let $\mathcal P'$ consist of the four triples $(\kappa,\gamma,b')$ in Table~\ref{tab:extra-root}. For the pairs in $\mathcal P$, let
\begin{align}
 B_\gamma(d)=b_\gamma-\gamma(1-d),\qquad
 B_*(d)=\min_{(\gamma,b_\gamma)\in\mathcal P} B_\gamma(d).
 \label{eq:B-gamma}
\end{align}

Let
\begin{align}
 \tau=\ell_1<\ell_2<\cdots<\ell_m<\ell_{m+1}=1,\qquad m=176,
 \label{eq:intervals}
\end{align}
be the degree cutoffs in Table~\ref{tab:cutoffs}. Let $F_i$ consist of the edges $e\in F$ with $\widehat d_e\in[\ell_i,\ell_{i+1})$, with the right endpoint included for $i=m$. The degree classes $F_1,\ldots,F_m$ are fixed throughout the coloring.

\subsection{Weighted neighborhood sizes and densities}

For a root $f\in F$, let
\begin{align}
 y_i=\frac{|N_Y(f)\cap F_i|}{D},\quad
 z_i=\frac{|U_2(f)\cap F_i|}{D},\quad
 d=\widehat d_f.
 \label{eq:yz}
\end{align}
Then $\sum_i y_i=d+O(D^{-1})$, $0\le z_i\le y_i$, and $\sum_i z_i\le1-d+O(D^{-1})$ by Lemma~\ref{lem:incidence}. Let $w_1\le\cdots\le w_m=1$ be numbers in $(0,1]$. Let $R$ be the scale at the current stage. In class $i$, the prescribed proportion of vertices to retain uncolored, relative to the start of the procedure, is $R\widetilde w_i$, where $0<\widetilde w_i\le w_i$ and $\widetilde w_m=1$. Let $w=w_1$ and let $u\in[0,1]$ satisfy
\begin{align}
 w+(1-w)u=\frac{\sum_i w_i y_i}{d_f}.
 \label{eq:u}
\end{align}
If $w=1$, let $u=0$. Let $T=d(w+(1-w)u)$.

For $\lambda\le0$, let
\begin{align}
 A_\lambda&=\max_i\{a-\ell_i-\lambda w_i\},\quad
 B_\lambda=\max\{0,-\lambda-A_\lambda\},\notag\\
 H_\lambda(d,T)&=\lambda T+d A_\lambda+(1-d)B_\lambda.
 \label{eq:AB}
\end{align}
Then
\begin{align}
 \frac{1}{D}\sum_{e\in U_1(f)}(a-\widehat d_e)
 \le H_\lambda(d,T)+o(1).
 \label{eq:linear-bound}
\end{align}
Indeed, the left-hand side is at most $\sum_i(a-\ell_i)(y_i-z_i)$. By \eqref{eq:AB}, this is at most
\begin{align*}
 \sum_i(A_\lambda+\lambda w_i)(y_i-z_i)
  =dA_\lambda+\lambda T
    +\sum_i(-A_\lambda-\lambda w_i)z_i+o(1).
\end{align*}
Since $\lambda\le0$ and $w_i\le1$, each coefficient in the last sum is at most $B_\lambda$. This proves \eqref{eq:linear-bound}.

For a finite set $\Lambda\subseteq(-\infty,0]$ containing $0$, let
\begin{align}
 C(d,T)=\min\bigl\{d^2,\ &
 B_\gamma(d)+kH_\lambda(d,T),\
 b'-\gamma'(1-d)+\kappa H_\lambda(d,T) \notag\\[-2pt]
 &\quad : (\gamma,b_\gamma)\in\mathcal P,\ 
 (\kappa,\gamma',b')\in\mathcal P',\ \lambda\in\Lambda\bigr\}.
 \label{eq:C}
\end{align}
Lemma~\ref{lem:density} gives $c_f\le C(d,T)+o(1)$. The second row of Table~\ref{tab:extra-root} gives
\begin{align}
 C(d,T)\le 0.72861074<0.75.
 \label{eq:C-global}
\end{align}

Let $e_{ij}$ be the number of ordered edges of $Y[N_Y(f)]$ with first endpoint in $F_i$ and second endpoint in $F_j$, divided by $D^2$. Thus $0\le e_{ij}\le y_i y_j$ and $\sum_{i,j}e_{ij}=c_f$. For $V,A,B\ge0$, let
\begin{align}
 h&=\min\{V,A\},\quad j=\min\{V-h,B\},\notag\\
 \mathcal A_w(V;A,B)&=h+w j+w^2(V-h-j).
 \label{eq:allocation}
\end{align}

\begin{lemma}\label{lem:weighted}
We have
\begin{align}
 \sum_{i,j}w_iw_j e_{ij}
 \le \mathcal A_w\bigl(c_f;d^2u^2,2d^2u(1-u)\bigr)+o(1).
 \label{eq:weighted-bound}
\end{align}
The weighted count is also at most $T^2+o(1)$.
\end{lemma}

\begin{proof}
For any $z\ge0$,
\begin{align*}
 \sum_{i,j}w_iw_j e_{ij}
 \le zc_f+\sum_{i,j}y_i y_j[w_iw_j-z]_+.
\end{align*}
Let $W,W'$ be independent random variables taking value $w_i$ with probability $\frac{y_i}{d_f}$. For fixed $W'$, the function $W\mapsto[WW'-z]_+$ is convex. Replacing $W$ by a variable on $\{w,1\}$ with the same mean increases the expectation. Apply the same argument to $W'$. The sum is therefore at most
\begin{align*}
 d_f^2\bigl(u^2[1-z]_+
       +2u(1-u)[w-z]_++(1-u)^2[w^2-z]_+\bigr).
\end{align*}
Minimizing over $z$ is equivalent to allocating $c_f$ to the three capacities $d_f^2u^2$, $2d_f^2u(1-u)$, and $d_f^2(1-u)^2$, with respective weights $1,w,w^2$. The maximum fills the capacities in this order. Since $c_f\le d_f^2$, the remaining mass fits in the third capacity and receives weight $w^2$. Replacing $d_f$ by $d$ proves \eqref{eq:weighted-bound}. The other bound follows from $e_{ij}\le y_i y_j$ and $\sum_i w_i y_i=T+O(D^{-1})$.
\end{proof}

Let $J$ be the graph induced by the uncolored vertices of $Y$ at a stage with scale $M=RD+o(D)$, and let $S=N_J(f)$ for an uncolored vertex $f$. Let $t,c$ be as in \eqref{eq:tc}, and let
\begin{align}
 m_f=\frac{1}{D}\sum_{e\in N_Y(f)}x_e,\qquad
 \mu=\frac{1}{M}\sum_{e\in S}x_e.
 \label{eq:m-and-mu}
\end{align}
The concentration estimates proved in Section~\ref{sec:coloring} give, uniformly at each stage,
\begin{align}
 t=\sum_i\widetilde w_i y_i+o(1),\quad
 c=\sum_{i,j}\widetilde w_i\widetilde w_j e_{ij}+o(1),\quad
 \mu=\frac{1}{D}\sum_i\widetilde w_i
                 \sum_{e\in N_Y(f)\cap F_i}x_e+o(1).
 \label{eq:retained-counts}
\end{align}
In particular, $t\le T+o(1)$ and $\mu\le m_f+o(1)$. For the upper bound on $c$, we use Lemma~\ref{lem:weighted} with the upper weights $w_i$. Let
\begin{align*}
 c^*(d,u)=\min\bigl\{T^2,\mathcal A_w([C(d,T)]_+;d^2u^2,2d^2u(1-u))\bigr\}.
\end{align*}
Then $c\le c^*(d,u)+o(1)$, and \eqref{eq:C-global} ensures that $1-c\ge0.25$ for sufficiently large $\Delta$.

\subsection{The sum of the degree excesses}

For $d\in[\tau,1]$, let
\begin{align}
 G_*(d)=B_*(d)+k(a-\tau)d+k(1-a)(1-d).
 \label{eq:Gstar}
\end{align}
Lemma~\ref{lem:density} gives
\begin{align}
 m_f\le\frac{G_*(d)-c_f}{k}+o(1).
 \label{eq:degree-excess}
\end{align}
To see this, let $z=\frac{|U_2(f)|}{D}$. The signed sum in \eqref{eq:linear-bound} equals
\begin{align*}
 (a-\tau)(d-z)-\frac{1}{D}\sum_{e\in U_1(f)}x_e+o(1).
\end{align*}
The contribution to $m_f$ from $U_2(f)$ is at most $(1-\tau)z$. Substitution in the bound with $B_*(d)$, followed by $z\le1-d+o(1)$, proves \eqref{eq:degree-excess}.

Consider a rectangle $[d_-,d_+]\times[u_-,u_+]$ in the parameters $(d,u)$. Let $t_+$, $V$, $A_+$, and $B_+$ be upper bounds throughout the rectangle on $T$, $[C(d,T)]_+$, $d^2u^2$, and $2d^2u(1-u)$, respectively. Let $C_+$ be an upper bound on $\mathcal A_w(V;A_+,B_+)$. By Lemma~\ref{lem:weighted} and \eqref{eq:retained-counts}, we have $c\le\mathcal A_w(c_f;A_+,B_+)+o(1)$. Since $c_f\le V+o(1)$ and $\mathcal A_w$ has slopes at most $1$ in its first argument, the same bound holds with $c_f$ replaced by $\min\{c_f,V\}$. The slopes are also at least $w^2$, so
\begin{align*}
 C_+-c
 &\ge \mathcal A_w(V;A_+,B_+)-\mathcal A_w(\min\{c_f,V\};A_+,B_+)-o(1)\\
 &\ge w^2(V-c_f)-o(1).
\end{align*}
Dividing by $w^2$ and rearranging gives
\begin{align}
 c_f\ge V-\frac{C_+-c}{w^2}-o(1).
 \label{eq:initial-density-lower}
\end{align}
We use $C_+$ before taking the minimum with $t_+^2$ in \eqref{eq:initial-density-lower}.

Let $G_+$ be an upper bound on $G_*(d)$ throughout the rectangle. It is obtained by taking the maximum of each affine expression in \eqref{eq:Gstar} at the two degree endpoints and then taking their minimum. Equations \eqref{eq:degree-excess} and \eqref{eq:initial-density-lower} give
\begin{align}
 \mu\le m_0+\frac{C_+-c}{kw^2}+o(1),\qquad
 m_0=\frac{[G_+-V]_+}{k}.
 \label{eq:mu-budget}
\end{align}

For each $(\kappa,\gamma,b')\in\mathcal P'$ with $\kappa>0$, let $G_+^{\kappa,\gamma,b'}$ be an upper bound throughout the rectangle on the affine function
\begin{align}
 b'-\gamma(1-d)+\kappa(a-\tau)d+\kappa(1-a)(1-d).
 \label{eq:extra-excess-affine}
\end{align}
The proof of \eqref{eq:degree-excess} applies to \eqref{eq:density-extra}. Thus \eqref{eq:initial-density-lower} gives
\begin{align}
 \mu\le\frac{[G_+^{\kappa,\gamma,b'}-V]_+}{\kappa}
       +\frac{C_+-c}{\kappa w^2}+o(1).
 \label{eq:extra-mu-budget}
\end{align}
We also have $\mu\le(1-\tau)t_++o(1)$. If $w<1$, let
\begin{align*}
 \xi=\max_{1\le i\le m}
 \left\{w_i(\ell_{i+1}-\tau)-\frac{(1-\tau)(w_i-w)}{1-w}\right\}.
\end{align*}
The term with $i=m$ is zero, so $\xi\ge0$. We have $\mu\le d_+(\xi+(1-\tau)u_+)+o(1)$. Indeed, \eqref{eq:retained-counts} and $x_e\le\ell_{i+1}-\tau$ for $e\in F_i$ give $\mu\le\sum_iw_i(\ell_{i+1}-\tau)y_i+o(1)$. The definition of $\xi$ and \eqref{eq:u} bound this sum by $d(\xi+(1-\tau)u)+o(1)$.

\subsection{Bounds on rooted triangle counts}

After applying Lemma~\ref{lem:split}, enlarge $S$ to $M$ candidates as in Section~\ref{sec:offers}. Let $H$ be the graph on the enlarged set, with the added vertices isolated, and let $Q=\overline H$. For a candidate $v$, let
\begin{align}
 r_v=\frac{d_H(v)}{M},\quad
 a_v=\frac{2e(H[N_H(v)])}{M^2},\quad
 b_v=\frac{6\,(\text{number of triangles in }H[N_H(v)])}{M^3}.
 \label{eq:rooted-moments}
\end{align}
The last two quantities count triangles and copies of $K_4$ through $v$, respectively, with the other vertices ordered. Let $\langle Z\rangle=\frac{1}{M}\sum_{v\in V(H)}Z_v$. For an added candidate, let $x_v=0$. Then $\langle x\rangle=\mu$.

We use bounds of the form
\begin{align}
 a_v\le U+\nu x_v+o(1),\qquad \nu\ge0.
 \label{eq:root-triangle-cap}
\end{align}
The pair $(U,\nu)$ is common to all candidates at the given upper weights. We obtain such pairs in three ways.

First, $a_v$ is at most the normalized ordered edge count in the whole current neighborhood of $v$. By \eqref{eq:retained-counts}, this is at most $c_v+o(1)$. Since $\widehat d_e\ge\tau$, Lemma~\ref{lem:density} gives \eqref{eq:root-triangle-cap} with
\begin{align}
 U_\gamma=B_\gamma(\tau)+k(a-\tau)\tau,\qquad
 \nu_\gamma=\gamma+k(a-\tau).
 \label{eq:initial-triangle-cap}
\end{align}

Second, let $\beta_w=\max_{i>1}\frac{w_i-w}{\ell_i-\tau}$. The weighted neighborhood density at $v$ is at most
\begin{align}
 w^2c_v+(1+w)\beta_w \widehat d_v m_v+o(1).
 \label{eq:weighted-triangle-intermediate}
\end{align}
Indeed, $w_iw_j-w^2\le(w_i-w)+w(w_j-w)$, and $w_i-w\le\beta_w x_e$ for $e\in F_i$. The row sums of the ordered edge counts at $v$ are at most $\widehat d_v y_i$, where $y_i$ is evaluated at $v$. By symmetry this gives \eqref{eq:weighted-triangle-intermediate}.

Let
\begin{align*}
 h_v=\frac{1}{D}\sum_{e\in U_1(v)}x_e
     +(a-\tau)\frac{|U_2(v)|}{D}.
\end{align*}
Lemma~\ref{lem:density} gives $c_v\le B_\gamma(\widehat d_v)+k(a-\tau)\widehat d_v-kh_v+o(1)$, and $m_v\le h_v+(1-a)(1-\widehat d_v)+o(1)$. If $(1+w)\beta_w\le kw^2$, the coefficient of $h_v$ after substitution in \eqref{eq:weighted-triangle-intermediate} is nonpositive. Also $\widehat d_v(1-\widehat d_v)\le\tau(1-\tau)$ since $\widehat d_v\ge\tau>\frac{1}{2}$. Thus \eqref{eq:root-triangle-cap} holds with
\begin{align}
 U=w^2U_\gamma+(1+w)\beta_w(1-a)\tau(1-\tau),\qquad
 \nu=w^2\nu_\gamma.
 \label{eq:retention-triangle-cap}
\end{align}

Third, for a chosen $\nu\ge0$, we may take
\begin{align}
 U\ge\sup_{\substack{\tau\le d\le1\\0\le u\le1}}
                \bigl(c^*(d,u)-\nu(d-\tau)\bigr).
 \label{eq:scan-triangle-cap}
\end{align}
This follows from \eqref{eq:retained-counts} and Lemma~\ref{lem:weighted} at $v$. The calculation uses $\nu\in\{0.48,0.55,0.60,0.65,0.70,0.74,0.80,0.88,0.98\}$. Each of \eqref{eq:initial-triangle-cap}, \eqref{eq:retention-triangle-cap}, and \eqref{eq:scan-triangle-cap} also gives the desired inequality for the added vertices, since $a_v=x_v=0$ there and the bounds $U$ are nonnegative.

For a true candidate $v$, we have $d_H(v)\le d_J(v)$, and \eqref{eq:retained-counts} at $v$ gives $\frac{d_J(v)}{M}\le d_v+o(1)$. For the added candidates, $r_v=x_v=0$. Hence for every candidate $v$,
\begin{align}
 r_v\le\tau+x_v+o(1).
 \label{eq:root-degree-cap}
\end{align}

\subsection{The offer probability with the triangle bounds}
\label{sec:triangle-offer}

In addition to $L_0$, we use the inequalities in Table~\ref{tab:moment-gains}. For each row, with the indicated weight $\omega_v$, constant $U_*$, and affine function $L_*$, Proposition~\ref{prop:ordinary-finite} gives
\begin{align}
 \mathcal F(Q)+\langle\omega_v(a_v-U_*)\rangle
 \ge L_*(s)-o(1),\qquad s=1-c,
 \label{eq:moment-affine}
\end{align}
on the stated interval. For consecutive rows with the same weight and $U_*$, let $L_*$ be the resulting continuous piecewise affine function. For the weight $0.24r_v^3+0.025b_v$, interpolation of the inequalities for $U_*=0.558$ and $U_*=0.618$ gives an inequality of the same form at each intermediate $U_*$.

Consider a rectangle as in \eqref{eq:mu-budget}, and let $(U,\nu)$ satisfy \eqref{eq:root-triangle-cap}. Let $\phi$ be an upper bound on $\omega_v$, let $W_0$ be an upper bound on $\langle\omega\rangle$, and let $A_0(c)$ be a lower bound on $\langle\omega\rangle$. For $\omega_v=0.1a_v$, we use
\begin{align}
 \phi=0.1 C_+,\qquad W_0=0.1,\qquad A_0(c)=0.1c(2c-1).
 \label{eq:triangle-weight-bounds}
\end{align}
To prove $A_0(c)\le\langle\omega\rangle$, sum $|N_H(u)\cap N_H(v)|\ge d_H(u)+d_H(v)-M$ over ordered adjacent pairs $u,v$. The sum on the left is $M^3\langle a\rangle$, and Cauchy-Schwarz bounds the sum of the degree terms below by $2M^3c^2$.

For $\omega_v=\eta_3r_v^3+\eta_4b_v$, we use
\begin{align}
 \phi=\eta_3t_+^3+\eta_4 C_+^{\frac{3}{2}},\qquad
 W_0=\eta_3+\eta_4,\qquad A_0(c)=\eta_3c^3.
 \label{eq:star-weight-bounds}
\end{align}
Here $r_v\le t_++o(1)$ and $a_v\le c$. The trace bound for triangles gives $b_v\le a_v^{\frac{3}{2}}$. Since $\langle r\rangle=c$, Jensen's inequality gives $\langle r^3\rangle\ge c^3$. In these estimates we absorb the uniform $o(1)$ errors in the upper bounds on $t$ and $c$.

Let $g_0=[U_*-U]_+$ and $h_0=[U-U_*]_+$. By \eqref{eq:root-triangle-cap} and \eqref{eq:mu-budget},
\begin{align}
 \langle\omega_v(a_v-U_*)\rangle
 &\le \nu\phi\mu+W_0h_0-g_0 A_0(c)+o(1)\notag\\
 &\le \nu\phi m_0+W_0h_0-g_0 A_0(c)
             +\frac{\nu\phi}{kw^2}(C_+-c)+o(1).
 \label{eq:moment-correction}
\end{align}
Let $[s_-,s_+]$ be the interval of the chosen function $L_*$, and let $\kappa$ be a lower bound on its slopes. Let $s_0=1-C_+$. Recall that $s=1-c$ and that $\frac{2e(Q)}{M^2}=s-\frac{1}{M}$. Assume $s_0\in[s_-,s_+]$ and
\begin{align}
 \kappa\ge\frac{\nu\phi}{kw^2}+g_0 A_0'(C_+).
 \label{eq:slope-test}
\end{align}
Since $c\le C_++o(1)$, we have $s\ge s_0-o(1)$. If $s\le s_+$, let $s'=\max\{s,s_-\}$. Then $|s'-s|=o(1)$, and the endpoint estimate in Appendix~\ref{app:ordinary} allows us to use $L_*(s')$ in \eqref{eq:moment-affine}. The slope bound and the fact that $s'\ge s_0-o(1)$ give $L_*(s')-L_*(s_0)\ge\kappa(C_+-c)-o(1)$. Also $A_0(C_+)-A_0(c)\le A_0'(C_+)(C_+-c)+o(1)$, by \eqref{eq:triangle-weight-bounds} or \eqref{eq:star-weight-bounds}. Combining \eqref{eq:moment-affine} and \eqref{eq:moment-correction} with \eqref{eq:slope-test} gives
\begin{align}
 \mathcal F(Q)\ge
 L_*(1-C_+)-\nu\phi m_0-W_0h_0+g_0A_0(C_+)-o(1).
 \label{eq:moment-lower}
\end{align}
If $s>s_+$, the probability of no offer in the padded model is at least $L_0(s_+)-o(1)$, since $L_0$ is nondecreasing. We may therefore use the minimum of the right-hand side of \eqref{eq:moment-lower} and $0.18711$ for the star weights, or $0.20835$ for the weight $0.1a_v$. These caps are bounded above by $L_0(0.42)$ and $L_0(0.48)$, respectively.

For the weight $0.24r_v^3+0.025b_v$, we use $U_*=\max\{0.558,\min\{U,0.618\}\}$ and the interpolated function $L_*$ when $U<0.618$. Its slopes are at least $0.263$. For each pair $(U,\nu)$ and each weight, apply \eqref{eq:moment-lower} when $1-C_+\in[s_-,s_+]$ and \eqref{eq:slope-test} holds.

We also use the inequalities
\begin{align}
 \mathcal F(Q)+\langle\omega_v(a_v-U_*)+\zeta_v(r_v-R_*)\rangle\ge\iota+\rho s-o(1)\qquad(s\in[s_-,s_+]),
 \label{eq:mixed-moment-affine}
\end{align}
certified in Proposition~\ref{prop:ordinary-finite}. Here $\omega_v$ and $\zeta_v$ are nonnegative linear combinations of the following monomials, respectively:
\begin{align}
 &1,\ r,\ r^2,\ r^3,\ a,\ ra,\ b,\notag\\
 &1,\ r,\ r^2,\ r^3,\ r^4,\ a,\ ra,\ r^2a,\ a^2,\ b,\ rb.
 \label{eq:mixed-monomials}
\end{align}
The variables in \eqref{eq:mixed-monomials} are evaluated at the same vertex $v$. The rooted degree and triangle bounds can both be used for a moment such as $\langle ra^2\rangle$. Indeed, this moment uses the root, one further vertex for $r$, and two pairs of further vertices for $a^2$.

For a weight $\omega$ or $\zeta$, let $P_\omega$ or $P_\zeta$ denote its pointwise upper bound, let $W_\omega$ or $W_\zeta$ denote its upper bound on the average, and let $A_\omega(c)$ or $A_\zeta(c)$ denote its lower bound on the average. These are obtained by summing the corresponding bounds on the monomials with their nonnegative coefficients. For the pointwise bounds, use $r_v\le t_++o(1)$, $a_v\le\min\{C_+,t_+^2\}+o(1)$, and $b_v\le a_v\sqrt{C_+}+o(1)$. If a monomial $r^ia^jb^k$ has $m'=i+2j+3k>0$, it is at most $r^{m'}$. Its average is thus at most $C_+t_+^{m'-1}+o(1)$ and is also at most its pointwise upper bound.

For the lower bounds, Jensen's inequality on the true candidates gives $\langle r^j\rangle\ge\frac{c^j}{t_+^{j-1}}-o(1)$ for $1\le j\le4$. Let
\begin{align}
 A(c)=\left[\frac{2c^2}{t_+}-t_+c\right]_+.
 \label{eq:triangle-average-support}
\end{align}
On summing the common-neighbor bound over ordered adjacent pairs of $H$, with $|S|$ in place of $M$, we obtain $\langle a\rangle\ge 2\langle r^2\rangle-tc\ge A(c)-o(1)$. Here we also use $\langle a\rangle\ge0$. Since $a_v\le r_v^2$, Jensen's inequality on the true candidates gives the lower bounds $\frac{A(c)^{3/2}}{\sqrt{t_+}}-o(1)$ for $\langle ra\rangle$ and $\frac{A(c)^2}{t_+}-o(1)$ for both $\langle r^2a\rangle$ and $\langle a^2\rangle$. The lower bounds for $\langle b\rangle$ and $\langle rb\rangle$ are zero. For the constant monomial all three bounds are $1$. The resulting lower bound functions are convex for $c\ge0$.

Let $g_a=[U_*-U]_+$, $h_a=[U-U_*]_+$, $g_r=[R_*-\tau]_+$, and $h_r=[\tau-R_*]_+$. Let $\Gamma=\nu P_\omega+P_\zeta$ and $A_*(c)=g_aA_\omega(c)+g_rA_\zeta(c)$. By \eqref{eq:root-triangle-cap} and \eqref{eq:root-degree-cap},
\begin{align}
 \langle\omega_v(a_v-U_*)+\zeta_v(r_v-R_*)\rangle
 \le \Gamma\mu+h_aW_\omega+h_rW_\zeta-A_*(c)+o(1).
 \label{eq:mixed-correction}
\end{align}
The bounds on $\mu$ above give several upper bounds of the form $\Gamma\mu\le p+z(C_+-c)+o(1)$. Each convex combination of these upper bounds is also valid. Since $A_*$ is convex, if $1-C_+\in[s_-,s_+]$ and $\rho\ge z+A_*'(C_+)$, then \eqref{eq:mixed-moment-affine} and \eqref{eq:mixed-correction} give
\begin{align}
 \mathcal F(Q)\ge \iota+\rho(1-C_+)+A_*(C_+)-h_aW_\omega-h_rW_\zeta-p-o(1)
 \label{eq:mixed-moment-lower}
\end{align}
for $s$ in the interval. Here $A_*'$ denotes the right derivative, and the endpoint estimate applies as in \eqref{eq:moment-lower}. For $s>s_+$, we again use $L_0(s_+)-o(1)$. The right endpoints $s_+=0.42,0.48,0.60$ of the new inequalities have caps $0.18711,0.20835,0.25144$, respectively, each at most $L_0(s_+)$.

Let $c_+$ be the upward-rounded value of $\min\{C_+,t_+^2\}$ on the grid in Appendix~\ref{app:palette}, so that $c\le c_++o(1)$. Let $L$ be the largest of the resulting capped lower bounds and $L_0(1-c_+)$. By \eqref{eq:padding} and \eqref{eq:three-term},
\begin{align}
 \Prb(I\cap S\ne\varnothing)
 \le \min\left\{
 t_+-g([t_+^2-c_+]_+),\ 1-Le^{1-t_+}
 \right\}+o(1).
 \label{eq:rectangle-offer}
\end{align}
For the first term, the integer bound in \eqref{eq:rectangle-t} is at most $M_0$, so $t_+\le1$ in real units. On $0\le t\le1$ and $c\ge0$, the function $t-g([t^2-c]_+)$ is nondecreasing in $c$, since $g$ is nondecreasing on $[0,1]$. When $t^2>c$, its derivative with respect to $t$ is $1-2t g'(t^2-c)\ge1-t\ge0$. When $t^2\le c$, the function is $t$. By continuity it is nondecreasing in $t$ throughout the domain.

For each class $i$, we bound the probability in \eqref{eq:rectangle-offer} over the rectangle $[\ell_i,\ell_{i+1}]\times[0,1]$ by subdivision. Appendix~\ref{app:palette} describes the calculation.

\section{The coloring procedure}\label{sec:coloring}

The degree classes $F_1,\ldots,F_m$ and the neighborhoods in $Y$ are fixed. We first give the concentration estimates for a single round. We then choose the list sizes to obtain the prescribed uncolored proportions in these classes. In a fixed finite sequence of rounds, each round uses working lists of size $O(D)$. The number of distinct color labels is also $O(D)$, since a common collection of $O(D)$ fresh colors suffices to supply the lists in any one round.

\subsection{One coloring round}

Let $J$ be the graph induced by the uncolored vertices at the start of a round. Let $M\ge\Delta(J)$. Each uncolored vertex has a \emph{working list}, which is a prescribed subset of the colors available at that vertex under the current partial coloring. For each color in at least one working list, choose an independent set by the priority rule in Section~\ref{sec:offers}, independently of the other colors and of previous rounds. A vertex is offered that color if it belongs to the independent set and the color is in its working list. Each vertex receiving at least one offer keeps one of its offered colors. Since each color is offered on an independent set, the resulting partial coloring is proper.

At each vertex that remains uncolored, delete every color offered to any of its neighbors in the round from its working list. Each color remaining in the list is therefore available. We then add fresh colors, or trim the list, to attain the working list size prescribed for the next round.

\begin{lemma}\label{lem:round}
Let $M=\Theta(D)$ and assume that there are $O(D)$ colors in the working lists. Let the working lists of the vertices in class $i$ have size $L_i=O(D)$, and let
\begin{align}
 \mu_i=\left(1-\frac{1}{M+1}\right)^{L_i}.
 \label{eq:mu}
\end{align}
There is an outcome of the round for which the following hold simultaneously at every vertex $f$ of $Y$.
\begin{enumerate}
 \item The number of uncolored neighbors of $f$ in class $i$ and the sum of their weights $x_e$ each equal their values before the round multiplied by $\mu_i$, with error $O(D^{\frac34})$.
 \item The number of ordered edges between classes $i,j$ in the uncolored neighborhood of $f$ equals its value before the round multiplied by $\mu_i\mu_j$, with error $O(D^{\frac74})$.
 \item Let $f$ be an uncolored vertex in class $i$. For each color in the working list of $f$, assume that, conditional on the state at the start of the round, the probability that at least one uncolored neighbor of $f$ receives an offer of that color is at most $p_i+o(1)$. If $f$ remains uncolored, the number of colors deleted from that list is at most $L_i(p_i+o(1))+O(D^{\frac34})$.
\end{enumerate}
The errors are uniform for any fixed finite sequence of rounds.
\end{lemma}

\begin{proof}
A vertex remains uncolored precisely when it receives no offer. Lemma~\ref{lem:priority} and independence between colors give \eqref{eq:mu}. For two distinct vertices, the probability that both receive a given color is $O(M^{-2})$, whether or not their lists contain that color. Thus the probability that both remain uncolored differs from the product of their marginal probabilities by $O(M^{-1})$, by multiplication over the $O(M)$ colors. This gives the expectations of the two counts in part (1) and the count in part (2), with an additional error $O(D)$ in the edge count. For part (3), for each color in the working list of $f$, consider the indicator that a neighbor with that color in its working list receives it as an offer. If $f$ remains uncolored, the sum of these indicators counts the deleted colors. Conditional on the state at the start of the round, the indicators are independent over the colors. For each color, the event counted by its indicator implies that its independent set meets the uncolored neighborhood of $f$, so the probability bound in Section~\ref{sec:profiles} bounds its expectation by $p_i+o(1)$.

For a fixed vertex $f$ of $Y$ and a fixed color, let $S=N_Y(f)\cap V(J)$ be the uncolored neighborhood of $f$ at the start of the round, and let $r=\frac{\log^2 D}{M}$. The probability of an offer at a vertex with priority greater than $r$ is at most $e^{-Mr}$. By a union bound, the probability of such an offer anywhere in $S$ is at most $|S|e^{-Mr}$. The number of vertices of $S$ with priority at most $r$ has mean $O(\log^2 D)$ and is a sum of independent indicators. The union bound and a Chernoff bound for this count imply that, for every fixed $A$ and all sufficiently large $D$, the probability of more than $\log^3 D$ offers in $S$ is less than $D^{-A}$.

In applying the bounded differences inequality, replace the offers for a color in $S$ by the empty collection if their number exceeds $\log^3 D$. Changing the variables for one color then changes the counts in part (1) by $O(\log^3 D)$, since $0\le x_e\le1$, and the count in part (2) by $O(D\log^3 D)$. There are $O(D)$ independent color coordinates. The probability of a deviation of order $D^{\frac34}$ in part (1) or $D^{\frac74}$ in part (2) is at most
\begin{align*}
 \exp\left(-\Omega\left(\frac{D^{\frac12}}{\log^6 D}\right)\right).
\end{align*}
The truncation probability and its contribution to the expectations are smaller than every fixed negative power of $D$. A Chernoff bound gives the required concentration in part (3).

Each bad event depends only on priorities at distance at most $2$ from its root in $Y$, together with the auxiliary comparison variables at the neighbors. Events whose roots are at distance greater than $4$ have disjoint underlying variables. Since there are a fixed number of events per root and $\Delta(Y)\le D$, the dependency degree is polynomial in $D$. The symmetric Lov\'asz local lemma gives an outcome avoiding all bad events.
\end{proof}

\subsection{Prescribed removal rates}

Let the target fraction of uncolored vertices in class $i$ be $R\widetilde w_i$, where $\widetilde w_i\le w_i$ and $w_m=1$. Start with $R=1$ and $\widetilde w_i=1$ for all $i$. For a target retention factor $\mu_i$ in a round, use a working list of size
\begin{align}
 L_i=\ceil{\frac{\log\mu_i}
                 {\log(1-\frac{1}{M+1})}}.
 \label{eq:target-list}
\end{align}
The resulting retention factor differs from its target by $O(M^{-1})$. Iteration of Lemma~\ref{lem:round} over a fixed number of rounds gives, for each $f$, the uncolored neighbor counts and sums of weights in each class, and the ordered edge counts between classes. More precisely, the initial counts and sums in class $i$ are multiplied by $R\widetilde w_i$ with error $O(D^{\frac34})$, and the initial ordered edge counts between classes $i,j$ are multiplied by $R^2\widetilde w_i\widetilde w_j$ with error $O(D^{\frac74})$. Let
\begin{align}
 M=\ceil{RD+D^{\frac9{10}}}.
 \label{eq:current-M}
\end{align}
Summing the neighbor counts over the fixed number of classes gives $d_J(f)\le R d_Y(f)+O(D^{\frac{3}{4}})\le RD+O(D^{\frac{3}{4}})$ for every $f\in V(J)$. Hence $M\ge\Delta(J)$ for sufficiently large $D$. Dividing by $M$ or $M^2$ gives \eqref{eq:retained-counts}. All scales $R$ in the argument are positive constants fixed before $D$ tends to infinity.

Consider a phase starting at $R=R_0$ with upper weights $w_1\le\cdots\le w_m=1$. Let $p_i$ be upper bounds on the quantities in \eqref{eq:rectangle-offer}, with the $o(1)$ terms omitted, for the respective classes at these weights. By replacing each bound by the maximum of the bounds up to its index, we may assume that $p_1\le\cdots\le p_m$. Let $p_*=p_m$ and $\alpha_i=\frac{p_*}{p_i}$. The calculation uses $1\le\alpha_i<2$. For an integer $q\ge2$, prescribe the relative retention
\begin{align}
 \widetilde w_i(R)=\widetilde w_i(R_0)
       \left(\frac{R}{R_0}\right)^{\alpha_i-1},
 \qquad \left(1-\frac{1}{q}\right)R_0\le R\le R_0.
 \label{eq:trajectory}
\end{align}
Each $\widetilde w_i$ is at most its initial upper bound throughout the phase. The bounds $p_i+o(1)$ therefore apply in every round of the phase by Section~\ref{sec:profiles}.

\begin{lemma}\label{lem:phase}
For every $\varepsilon>0$, the phase in \eqref{eq:trajectory} can be implemented with at most
\begin{align}
 \frac{p_*R_0D}{q}+\varepsilon D+o(D)
 \label{eq:phase-cost}
\end{align}
new colors.
\end{lemma}

\begin{proof}
Divide the interval in \eqref{eq:trajectory} into finitely many steps of equal logarithmic length $h$. In each step, the target retention factor of class $i$ is $e^{-\alpha_i h}$ because $R\widetilde w_i(R)$ is proportional to $R^{\alpha_i}$. Hence \eqref{eq:target-list} gives $L_i=\alpha_i hM+O(1)$. Successive working list sizes in the phase differ by $O(h^2D)+o(D)$. By Lemma~\ref{lem:round}, at most $L_i(p_i+o(1))+O(D^{\frac34})$ colors are deleted from each working list in class $i$. A common collection of at most
\begin{align*}
 hM\max_i(\alpha_i p_i)+O(h^2D)+o(D)
   =hMp_*+O(h^2D)+o(D)
\end{align*}
fresh colors restores every prescribed list size after trimming any excess. The initial working lists of the phase use $O(hD)+o(D)$ colors.

Summing over the steps and choosing $h$ sufficiently small bounds the cost by
\begin{align*}
 p_*D\int_{(1-\frac1q)R_0}^{R_0} dR+\varepsilon D+o(D).
\end{align*}
The number of steps is fixed before $D$ tends to infinity, so the concentration estimates apply in every step. This proves \eqref{eq:phase-cost}.
\end{proof}

\begin{lemma}\label{lem:tail}
At a stage of scale $RD$ with upper weights $w_i$, let $p_0$ be a common upper bound for all classes obtained from \eqref{eq:rectangle-offer} at these weights, with the $o(1)$ term omitted. For every $\varepsilon>0$, the partial coloring can be completed with at most $p_0 RD+\varepsilon D+o(D)$ new colors.
\end{lemma}

\begin{proof}
Use the same retention factor in every class, so the relative retentions remain constant. The same upper weights $w_i$ apply in every round, so \eqref{eq:rectangle-offer} continues to give the bound $p_0+o(1)$. The argument in Lemma~\ref{lem:phase}, with all logarithmic rates equal to $1$, takes the scale from $RD$ to $rD$ using at most $p_0(R-r)D+\frac{\varepsilon}{2}D+o(D)$ new colors. A greedy coloring of the remaining graph uses at most $rD+o(D)+1$ further colors. Choose $0<r<\min\{R,\frac{\varepsilon}{2}\}$ and then let $D$ tend to infinity.
\end{proof}

\subsection{The palette bound}

Appendix~\ref{app:palette} specifies the upper weights, the bounds $p_i$ for each phase, and the number of phases. Let $R_j$ be the starting scale of phase $j$, let $q_j$ be its phase parameter, and let $p_{*,j}$ be its largest offer probability bound. Let $R_{\mathrm{end}}$ and $p_0$ be the scale and the common probability bound at the completion step, respectively.

\begin{proposition}\label{prop:palette}
The recurrence in Appendix~\ref{app:palette} gives $487$ phases with $0<p_i\le p_{*,j}<2p_i$ in each phase. The normalized leading terms of the phase and completion costs satisfy
\begin{align}
 \sum_{j=0}^{486}\frac{p_{*,j}R_j}{q_j}
   +p_0 R_{\mathrm{end}}
 \le \frac{819392866996}{10^{12}}.
 \label{eq:palette-total}
\end{align}
\end{proposition}

\begin{proof}
The class bounds follow by applying \eqref{eq:rectangle-offer} on the rectangles in Appendix~\ref{app:palette}. Every upper weight is rounded up along the trajectory \eqref{eq:trajectory} at the end of each phase. The recurrence rounds upper bounds on the normalized leading costs of the phases and the completion step up to multiples of $\frac{1}{10^{12}}$ and adds them. The exact comparison \eqref{eq:integer-output} gives the claimed value.
\end{proof}

\section{Proof of Theorem~\ref{thm:main}}\label{sec:mainproof}

\begin{proof}[Proof of Theorem~\ref{thm:main}]
By Lemma~\ref{lem:regular}, we may assume that $G$ is $\Delta$-regular. Let $D=2\Delta^2$, let $\tau=0.8175$ and $\tau_1=0.8194$, and let $K=\floor{\tau_1D}$. Repeatedly delete vertices of current degree less than $K$ from $X=L(G)^2$, and let $Y=X[F]$ be the remaining graph. If $F=\varnothing$, then Lemma~\ref{lem:deletion} gives a proper coloring of $X$ with at most $K$ colors, as required.

Assume $F\ne\varnothing$. Every vertex of $Y$ has degree at least $K\ge\tau_1D-1\ge\tau D-1$, so the bounds in Lemma~\ref{lem:density} apply. Use the degree classes in \eqref{eq:intervals}, and apply Lemmas~\ref{lem:phase} and \ref{lem:tail} with the parameters of Proposition~\ref{prop:palette}. By \eqref{eq:palette-total},
\begin{align*}
 \frac{819392866996}{10^{12}}<0.8194.
\end{align*}
Choose the error allowances in the finitely many applications of the two lemmas so that their sum is less than half the gap in this inequality. For sufficiently large $\Delta$, the uniform $o(D)$ terms and the difference between $\tau_1D$ and $K$ are also absorbed by the gap. Thus $Y$ has a proper coloring with at most $K$ colors. By Lemma~\ref{lem:deletion}, the coloring extends to $X$ with at most $K$ colors. Therefore
\begin{align*}
 \chi'_s(G)=\chi(X)\le K\le 1.6388\Delta^2.
\end{align*}
\end{proof}

\section*{Declaration of generative AI use}
During the preparation of this work, the author used OpenAI's GPT-6.2 Sol and GPT-6.1 Sol to explore potential approaches, assist with calculations, and refine the clarity of the exposition. The author independently reviewed, verified, and edited all AI-generated content and takes full responsibility for the final manuscript.

\appendix
\section{Local flag calculations}\label{app:local}

This appendix gives the flag expressions and coefficient tests used in Lemma~\ref{lem:density}. The witnesses and the coefficient data for the calculations are supplied in \cite{certificates}.

\subsection{Flags, products, and averaging}

A colored graph has vertex colors $0,1$, meaning black and red, and edge states $0,1,2$, meaning absent, marked, and unmarked. An isomorphism preserves all colors and states. Let $\mathcal H_n$ be the isomorphism classes of colored graphs on $n$ vertices for which each component contains a black vertex. Let $\mathcal H_n^+$ be the subclass where every present red-red edge is marked. Direct enumeration gives
\begin{align}
 |\mathcal H_5|=17950,\qquad |\mathcal H_5^+|=9530.
 \label{eq:basis-count}
\end{align}
The enumeration is over $2^n3^{\binom n2}$ labeled colored graphs, with the locality test followed by identification under permutations of the vertices.

A \emph{type} $\sigma$ is a colored graph whose vertices are labeled $1,\ldots,t$. A \emph{$\sigma$-flag} consists of a colored graph together with a distinguished induced copy of $\sigma$, whose labels are fixed by isomorphisms. Let $\mathcal H_n^\sigma$ be the local $\sigma$-flags on $n$ vertices. For a fixed induced labeled copy of $\sigma$, let the density of a flag with $r$ unlabeled vertices be its number of induced extensions divided by $\binom{\Delta}{r}$. We extend density evaluations linearly to expressions in flags. The notation $\varphi$ denotes such an evaluation on the colored graph under consideration, with the copy of the type fixed when labels are present.

The product of two $\sigma$-flags $F_1,F_2$ with respectively $r_1,r_2$ unlabeled vertices is determined by the coefficient
\begin{align}
 p(F_1,F_2;H)=\frac{1}{\binom{r_1+r_2}{r_1}}
 \sum_{\substack{A\dot\cup B=V(H)\setminus V(\sigma)\\|A|=r_1}}
 \ind_{H[V(\sigma)\cup A]\cong F_1}
 \ind_{H[V(\sigma)\cup B]\cong F_2}.
 \label{eq:local-product}
\end{align}
Products of linear combinations are obtained by bilinearity. The averaging operator $\av{\cdot}$ removes the labels. For a $\sigma$-flag $F$ on $n$ vertices, its coefficient on the underlying unlabeled graph is the proportion of injections of the $t$ labels into that graph that induce $F$. The denominator is $\frac{n!}{(n-t)!}$. These conventions agree with \cite{local,strong}.

For each fixed flag, the discrepancy between products of densities and densities of products is $O(\Delta^{-1})$. Terms with a repeated additional vertex have one fewer free choice. Locality and $|B|\le2\Delta$ bound their number by the required power of $\Delta$. The averaging of a square is asymptotically nonnegative. It can be verified directly by summing the square of a linear combination of conditional extension densities over the embeddings of the type. Allowing repeated extension vertices makes this a sum of squares, and restricting to distinct extension vertices changes it by $O(\Delta^{-1})$.

\subsection{The objective and the two degree expressions}

Let $H\in\mathcal H_5$, with vertex set $V$. Let $b_x$ be the indicator that $x$ is black, let $a_{xy}$ be the indicator that $xy$ is present, and let $f_{xy}$ be the indicator that $xy$ is marked. For disjoint unordered pairs $xy,uv$, let
\begin{align*}
 q_{xy,uv}=\ind_{\{xu,xv,yu,yv\}\cap E(H)\ne\varnothing}.
\end{align*}
Every sum in this subsection uses distinct vertices. For an ordered pair $(x,y)$ and an unordered pair $\{u,v\}$ in $V\setminus\{x,y\}$, let $w$ be the remaining vertex. Let
\begin{align}
 o(H)=\sum_{(x,y),\{u,v\}}
 b_yf_{xy}f_{uv}a_{xw}q_{xy,uv}(2-b_x)
                         \ind_{b_u+b_v\ge1}.
 \label{eq:o-coeff}
\end{align}
Let
\begin{align}
 p(H)&=\sum_{(x,y),\{u,v\}}
           (1-b_x)b_yf_{xy}f_{uv}a_{xw}q_{xy,uv},\notag\\
 q(H)&=6\sum_{(x,y)}(1-b_x)b_yf_{xy}
              \prod_{z\in V\setminus\{x,y\}}a_{xz},\notag\\
 s(H)&=10000p(H)-17297q(H),\label{eq:s-coeff}\\
 v(H)&=\sum_{(x,y)}b_yf_{xy}(2-b_x)
              \prod_{z\in V\setminus\{x,y\}}a_{xz}.
 \label{eq:m-coeff}
\end{align}
The required size-$5$ expressions are
\begin{align}
 O&=\sum_{H\in\mathcal H_5}\frac{o(H)}{60}H,
                  \label{eq:flag-objective}\\
 S&=\sum_{H\in\mathcal H_5}\frac{s(H)}{5000}H,
                  \label{eq:flag-signed}\\
 m&=\sum_{H\in\mathcal H_5}\frac{v(H)}{80}H.
                  \label{eq:flag-mass}
\end{align}
These formulas also apply after restriction to $\mathcal H_5^+$.

For \eqref{eq:flag-mass}, a marked edge with a black endpoint contributes $2\binom{\Delta-1}{3}$ to the count in \eqref{eq:m-coeff}. A marked black-black edge has two orientations of weight $1$, and a marked red-black edge has one orientation of weight $2$. Division by $80\binom{\Delta}{5}$ gives the marked-edge count divided by $2\Delta^2$, with an error $O(\Delta^{-1})$. This proves \eqref{eq:root-mass} in the flag normalization. The counts in \eqref{eq:o-coeff} and \eqref{eq:s-coeff} similarly give \eqref{eq:bridge-o} and \eqref{eq:bridge-s}.

\subsection{Square matrices}

We use the black one-vertex type and one labeled representative of each of the $38$ isomorphism classes of local colored graphs of order $3$. For a type $\sigma$ of size $t\in\{1,3\}$, let $r=\frac{5+t}{2}$ and let $\mathbf f_\sigma$ be the vector of flags in $\mathcal H_r^\sigma$. The one-vertex type has $55$ flags in the corresponding vector. Each three-vertex type has $53$ flags in the corresponding vector. Let $M_\sigma(H)$ be the rational symmetric matrix determined by
\begin{align}
 \av{\mathbf f_\sigma\mathbf f_\sigma^{\mathsf T}}
        =\sum_{H\in\mathcal H_5}M_\sigma(H)H.
 \label{eq:local-M}
\end{align}
Equations~\eqref{eq:local-product} and the averaging rule give every entry by a finite count. For any positive semidefinite matrix $Q_\sigma$,
\begin{align}
 \sum_H\tr(Q_\sigma M_\sigma(H))\varphi(H)\ge-o(1).
 \label{eq:local-square}
\end{align}
The same statement holds on $\mathcal H_5^+$ by assigning density zero to the excluded flags.

\subsection{Regularity and black-vertex constraints}

For a local type $\sigma$ and a label $i$, let $e_i^\sigma$ be the sum of extensions by one vertex adjacent to label $i$. Let $b_1^\sigma$ be the sum of extensions by one black vertex. Regularity gives $\varphi(e_i^\sigma)=1+O(\Delta^{-1})$, and the size of $B$ gives $\varphi(b_1^\sigma)\le2$. We use the following forms of expressions.
\begin{enumerate}
\item For a local type $\sigma$ of size $4$ and labels $i,j$, use
\begin{align}
 \av{e_i^\sigma-e_j^\sigma}.\label{eq:regular-row}
\end{align}
Its limiting value is zero. Both directions are therefore permissible as nonnegative expressions.
\item For a local type $\sigma$ of size $4$, use
\begin{align}
 \av{2e_1^\sigma-b_1^\sigma}.\label{eq:black-row}
\end{align}
\item Let $\sigma_B$ be the black one-vertex type. Include
\begin{align}
 2-\av{(e_1^{\sigma_B})^4}.\label{eq:black-normalization}
\end{align}
For $j=2,3,4$, let $b_j^{\sigma_B}$ be the sum of extensions whose $j$ new vertices are all black. Include both signs of
\begin{align}
 \av{\bigl(b_j^{\sigma_B}-(b_1^{\sigma_B})^j\bigr)
                       (e_1^{\sigma_B})^{4-j}}.
 \label{eq:black-moments}
\end{align}
\end{enumerate}
Every expression has a size-$5$ expansion, with a constant term where one occurs explicitly. The black-vertex moment expressions have limiting value zero by disjoint sampling. The normalization in \eqref{eq:black-normalization} evaluates as $2-\frac{|B|}{\Delta}+o(1)$. Thus every expression in this family has nonnegative limiting value for a regular graph with at most $2\Delta$ black vertices. Let $\mathcal R_0$ be this family, together with the basis nonnegativity expressions $H$ for $H\in\mathcal H_5$.

\subsection{Additional constraints at minimum conflict degree}

Let $\theta=0.8175$ as in Section~\ref{sec:density}. On $\mathcal H_5^+$, we use the following additional expressions, whose evaluations are nonnegative up to $o(1)$ when every marked edge with a black endpoint has conflict degree at least $\theta D-1$.

For a local three-vertex pattern $\sigma$ on an ordered marked edge $xy$ with $y$ black and a third vertex $w$, require that $w$ is black or adjacent to $x$ or $y$. For $H\in\mathcal H_5^+$, let
\begin{align}
 R_\sigma(H)=
 \sum_{(x,y,w)\text{ inducing }\sigma}
 \left(2f_{uv}q_{xy,uv}-8\theta a_{xw}a_{xu}a_{xv}\right),
 \label{eq:degree-row}
\end{align}
where $\{u,v\}$ are the remaining vertices. All orderings of the first three vertices are included. The corresponding expression is asymptotically nonnegative. For fixed $x,y,w$, the first count is at least $2\theta D-O(\Delta)$, whereas the subtracted count is at most $8\theta\binom{\Delta}{2}$. If $w$ is not adjacent to $x$, the subtracted count is zero. Summing over the $O(\Delta^3)$ choices of the fixed local triple gives an error $O(\Delta^{-1})$ after normalization. Using full ordered patterns permits the merging of any patterns related by an automorphism of the marked edge.

Next include the affine expressions $m-\theta$ and $1-m$. For a local unlabeled three-vertex graph $K$, let
\begin{align}
 P_K(H)&=\sum_{\substack{A\subseteq V(H)\\|A|=3}}
  \ind_{H[A]\cong K}f_{uv}\ind_{b_u+b_v\ge1},\notag\\
 Q_K(H)&=\sum_{\substack{A\subseteq V(H)\\|A|=3}}
  \ind_{H[A]\cong K}\sum_{w\in A}a_{wu}a_{wv},
 \label{eq:product-counts}
\end{align}
where $\{u,v\}=V(H)\setminus A$. Include the expressions with coefficients
\begin{align}
 \frac{P_K(H)}{40}-\theta\frac{Q_K(H)}{30},\qquad
 \frac{Q_K(H)}{30}-\frac{P_K(H)}{40}.
 \label{eq:mass-products}
\end{align}
Their evaluations are respectively $(m-\theta)\varphi(K)+o(1)$ and $(1-m)\varphi(K)+o(1)$. The factors $40$ and $30$ follow by comparing $\binom{\Delta}{5}$ with $\binom{\Delta}{3}$ and the counts of marked edges and pairs of neighbors of a fixed vertex.

Finally, let $\sigma$ be a local four-vertex pattern on $x,y,u,v$, with $xy$ marked and $y$ black. Include the expression whose coefficient on $H$ is
\begin{align}
 T_\sigma(H)=
 \sum_{(x,y,u,v)\text{ inducing }\sigma}
      a_{xw}\bigl(2(1-\theta)-a_{yw}\bigr),
 \label{eq:common-row}
\end{align}
where $w$ is the remaining vertex. For each fixed four-vertex pattern, its extension count is nonnegative up to $O(1)$ by \eqref{eq:codegree} and regularity. There are $O(\Delta^4)$ choices of the pattern, so the normalized error is $O(\Delta^{-1})$. The minimum conflict degree also gives the following four families of nonnegative expressions. As above, nonnegativity is up to $o(1)$ after summing over the conditioning pattern and dividing by $\binom{\Delta}{5}$.

\begin{enumerate}
\item For the same three-vertex pattern as in \eqref{eq:degree-row}, include the row with coefficient
\begin{align}
 \sum_{(x,y,w)\text{ inducing }\sigma}
       (2f_{uv}q_{xy,uv}-8\theta a_{xu}a_{xv}).
 \label{eq:unrestricted-degree-row}
\end{align}
For each fixed triple, the first count is at least $2\theta D-O(\Delta)$, and the subtracted count is at most $8\theta\binom{\Delta}{2}$. In this row the two neighbors of $x$ are counted for every triple of the given pattern.

\item For each such triple, let $C=N_G(x)\cup N_G(y)$ in $G$. The number of marked edges with an endpoint in $C$ is at least $\theta D-1$. Subtracting it from $\Delta|C|-e(G[C])$ gives
\begin{align}
 \Delta |N_G(x)\cap N_G(y)| + e(G[C])
 &+|\{e\in E(G)\setminus F^+ : e\text{ has an endpoint in }C\}|
 \notag\\
 &\le 2(1-\theta)\Delta^2+O(\Delta).
 \label{eq:degree-deficit}
\end{align}
For $z=u,v$, let $c_z=\ind_{\{a_{xz}+a_{yz}>0\}}$. Include the row obtained by summing
\begin{align}
 4(1-\theta)a_{xu}a_{xv}
 &-a_{xu}a_{xv}(a_{yu}+a_{yv})-a_{uv}c_uc_v \notag\\
 &-(a_{uv}-f_{uv})\ind_{\{c_u+c_v>0\}}
 \label{eq:deficit-coeff}
\end{align}
over the triples of the given pattern. The second term counts a common neighbor of $x,y$ together with a neighbor of $x$. The third and fourth terms count the last two terms on the left of \eqref{eq:degree-deficit}. The first term sums to $2(1-\theta)\Delta^2+O(\Delta)$ for a fixed triple, so the row is nonnegative up to the stated error.

\item The bound $m\ge\theta-O(D^{-1})$ and the incidence bound $\Delta|B|$ give $|B|\ge2\theta\Delta-O(1)$. For each local unlabeled graph $K$ of order $4$, include
\begin{align}
 \sum_{\substack{A\subseteq V(H)\\ |A|=4,\ H[A]\cong K}}
 \left(4b_v-2\theta\sum_{x\in A}a_{xv}\right),
 \label{eq:black-lower-row}
\end{align}
where $v$ is the remaining vertex. The sum over extensions of a fixed copy of $K$ is $4|B|-8\theta\Delta+O(1)$.

\item For each local unlabeled graph $K$ of order $3$, include
\begin{align}
 \sum_{\substack{A\subseteq V(H)\\ |A|=3,\ H[A]\cong K}}
 \left(b_u\sum_{x\in A}a_{xv}+b_v\sum_{x\in A}a_{xu}
          -3(b_u+b_v)f_{uv}\right),
 \label{eq:incidence-upper-row}
\end{align}
where $\{u,v\}=V(H)\setminus A$. For a fixed copy of $K$, the positive terms count $3|B|\Delta+O(\Delta)$ extensions. The subtracted terms count three times the number of marked incidences at $B$, which is at most $3|B|\Delta$. In particular, a marked edge with two black endpoints contributes two incidences.
\end{enumerate}

Each counted extension is local by its adjacency to the conditioning pattern, its adjacency to a black vertex, or its own black vertex. The number of copies of any conditioning pattern on $r$ vertices is $O(\Delta^r)$ for $r=3,4$. Excluding repeated vertices and the vertices of the pattern therefore gives a normalized error $O(\Delta^{-1})$ in the four families. Let $\mathcal R_\theta$ consist of the restrictions of $\mathcal R_0$ to $\mathcal H_5^+$, the rows \eqref{eq:degree-row}, \eqref{eq:mass-products}, \eqref{eq:common-row}, \eqref{eq:unrestricted-degree-row}, \eqref{eq:deficit-coeff}, \eqref{eq:black-lower-row}, and \eqref{eq:incidence-upper-row}, and the two affine constraints on $m$.

\subsection{Rational feasibility and the residual bound}

For any expression $R\in\mathcal R_\theta$, let
\begin{align*}
 R=R_0+\sum_H R_H H
\end{align*}
be its constant and size-$5$ coefficients. Here $R_0$ is the constant coefficient. For $\kappa,\gamma\ge0$, let
\begin{align}
 A_{\kappa,\gamma}(H)=\frac{o(H)}{60}
              +\frac{\kappa}{60}\frac{s(H)}{5000}
              -16\gamma\frac{v(H)}{80}.
 \label{eq:A-gamma}
\end{align}
The spanning-tree coefficient can also be computed by the matrix-tree theorem. With $A_H=(a_{uv})$ and underlying degrees $d_H(v)$, it is
\begin{align*}
 h(H)=\det\bigl(\diag(d_H(v)+b_v)-A_H\bigr).
\end{align*}
A rational certificate consists of nonnegative numbers $\lambda_R$ and positive semidefinite matrices $Q_\sigma$. Given a rational number $\epsilon\ge0$, require
\begin{align}
 A_{\kappa,\gamma}(H)+\sum_R\lambda_RR_H
      +\sum_\sigma\tr(Q_\sigma M_\sigma(H))
       \le\epsilon h(H)\qquad(H\in\mathcal H_5^+).
 \label{eq:rational-certificate}
\end{align}
Let $C=\sum_R\lambda_RR_0$. Multiplication by $\varphi(H)$, summation, and Lemma~\ref{lem:tree} give
\begin{align}
 O+\frac{\kappa}{60}S-16\gamma m\le C+4802\epsilon+o(1).
 \label{eq:weak-duality}
\end{align}
For the signed bound without a minimum degree assumption, let $\kappa=k$ and $\gamma=0$, let the basis be $\mathcal H_5$, and use only $\mathcal R_0$. The same criterion applies.

\begin{proposition}\label{prop:local-finite}
For $\kappa=k$ and $\gamma=0$, the rational system \eqref{eq:rational-certificate} on $\mathcal H_5$ with $\mathcal R_0$ has a solution with $C+4802\epsilon<10.6434$. For $\kappa=k$ and each $(\gamma,b_\gamma)$ in Table~\ref{tab:root}, the system on $\mathcal H_5^+$ with $\mathcal R_\theta$ has a solution with $C+4802\epsilon<16(b_\gamma-\gamma)$. For each $(\kappa,\gamma,b')$ in Table~\ref{tab:extra-root}, the same system has a solution with $C+4802\epsilon<16(b'-\gamma)$. The matrices in these solutions are positive definite.
\end{proposition}

\begin{proof}
The coefficient data and rational witnesses are supplied in \cite{certificates}. The verifier reconstructs the flag basis and checks that the archived square blocks and linear rows of \cite{strong} agree with the expressions above, up to positive scaling and permutation of the bases. It obtains the objective and the additional constraint coefficients from the counts in this appendix.

The verification checks nonnegativity of the scalar multipliers and positive definiteness of each of the $39$ matrices. It then evaluates the left-hand side of \eqref{eq:rational-certificate} at each flag. In the first system, the constant term and an upper bound on the residual coefficient are, respectively,
\begin{align*}
 \frac{26607971199}{2500000000}
 \qquad\text{and}\qquad
 \frac{59}{40000000000}.
\end{align*}
Since $h(H)\ge1$, the latter also bounds the residual coefficient divided by $h(H)$. For the other systems, the comparison divides each positive residual coefficient by its $h(H)$ and takes their maximum, with $0$ included. The resulting $\epsilon$ and the constant term satisfy $C+4802\epsilon<16(b_\gamma-\gamma)$ for each row of Table~\ref{tab:root} and $C+4802\epsilon<16(b'-\gamma)$ for each row of Table~\ref{tab:extra-root}. These are rational comparisons on the $9530$ flags in $\mathcal H_5^+$.

For the matrix check, clear denominators to obtain an integer symmetric matrix $A$. Index its rows and columns by $0,\ldots,n-1$, let $A^{(0)}=A$, and let $p_{-1}=1$. At step $r$, let $p_r=A^{(r)}_{rr}$ and, for $i,j>r$, let
\begin{align}
 A^{(r+1)}_{ij}
  =\frac{p_rA^{(r)}_{ij}-A^{(r)}_{ir}A^{(r)}_{rj}}{p_{r-1}}.
 \label{eq:bareiss}
\end{align}
The divisions are exact when the preceding pivots are nonzero. The successive pivots are the leading principal minors of $A$, so their positivity proves positive definiteness by Sylvester's criterion. Thus the verification of all the witness conditions uses integer arithmetic after denominators have been cleared.
\end{proof}

\section{Six-vertex flag calculations}\label{app:ordinary}

We use ordinary flag products and the positivity of squares \cite{razborov} to bound the expressions in \eqref{eq:no-offer-F}, \eqref{eq:moment-affine}, and \eqref{eq:mixed-moment-affine}. For an unlabeled graph $H$ on $r$ vertices and a graph $Q$ on $n\ge r$ vertices, let $p(H,Q)$ be the number of $r$-element subsets of $V(Q)$ that induce a copy of $H$, divided by $\binom{n}{r}$. All products below can be interpreted by independent sampling, with an error $O(n^{-1})$ when repeated samples are excluded.

\subsection{The coefficient of the objective}

Let $H$ be a labeled graph on $\{1,\ldots,6\}$, and let $a_{ij}$ be its adjacency indicators. Let
\begin{align}
 A_2=a_{12},\quad
 A_3=a_{12}a_{13}a_{23},\quad
 A_4=a_{12}a_{13}a_{23}a_{14}a_{24}a_{34}.
 \label{eq:clique-indicators}
\end{align}
For the indicated indices, let
\begin{align}
 B_j&=(1-a_{1j})(1-a_{2j}),\quad 3\le j\le6,\notag\\
 Y_j&=\max\{0,2-a_{1j}-a_{2j}-a_{3j}\},\quad 4\le j\le6,\notag\\
 U_j&=\binom{4-\sum_{i=1}^4a_{ij}}2,\quad 5\le j\le6.
 \label{eq:ordinary-variables}
\end{align}
Let $c_{rs}=1$ when $r,s\ge0$ and $r+s<3$, and let
\begin{align*}
 c_{30}=6,\qquad c_{21}=c_{12}=c_{03}=3.
\end{align*}
Let empty products have value $1$, and let
\begin{align}
 f_2(H)&=A_2\sum_{r=0}^4\frac{1}{2^{r+1}}
                      \prod_{j=3}^{r+2}B_j,\notag\\
 f_3(H)&=-\frac{A_3}{6}
   \sum_{\substack{r,s\ge0\\r+s\le3}}
       \frac{c_{rs}}{2^r3^s}
       \left(\prod_{j=4}^{r+3}B_j\right)
       \left(\prod_{j=r+4}^{r+s+3}Y_j\right),\notag\\
 f_4(H)&=A_4\left(\frac{1}{30}
               +\frac{h_1}{180}U_5
               +\frac{h_2}{1080}U_5U_6\right).
 \label{eq:ordinary-polynomial}
\end{align}
Here $h_1=\frac{23}{12}$ and $h_2=\frac{355}{144}$. Let
\begin{align}
 f(H)=\frac{1}{6!}\sum_{\pi\in S_6}
       \bigl(f_2(H^\pi)+f_3(H^\pi)+f_4(H^\pi)\bigr),
 \label{eq:ordinary-objective}
\end{align}
where $H^\pi$ has vertex $i$ corresponding to $\pi(i)$ in $H$. Then
\begin{align}
 \mathcal F(Q)=\sum_{H\in\mathcal G_6}f(H)p(H,Q)+O(n^{-1}),
 \label{eq:ordinary-expansion}
\end{align}
where $\mathcal G_6$ is the family of the $156$ unlabeled graphs on six vertices.

To verify the coefficients in \eqref{eq:ordinary-polynomial}, each factor $B_j$ has expectation $b_{12}$ after conditioning on vertices $1,2$. Each factor $Y_j$ has conditional expectation $b_{12}+b_{13}+b_{23}-b_{123}$ after conditioning on the first three vertices. Each factor $U_j$ has conditional expectation $6\overline b$ after conditioning on a $4$-clique. Independent additional samples give products of these expectations. This proves \eqref{eq:ordinary-expansion} directly from \eqref{eq:F}.

\subsection{The square matrices and interval constraints}

For each isomorphism class of graphs of order $t=0,2,4$, take one labeled representative as a type and use the corresponding flags of order $r=\frac{6+t}{2}$. There are $1$, $2$, and $11$ classes of the three orders, respectively, giving $14$ square blocks. The corresponding matrix dimensions are $4$, $20$, $20$, and eleven copies of $16$. For each type $\sigma$, let $\mathcal G_r^\sigma$ be the flags with fixed labels on their type. For $F,G$ in this basis, let $M_\sigma(H)_{F,G}$ be
\begin{align}
 \frac{1}{6!}\sum_{\pi\in S_6}\ind_{H[\pi(1),\ldots,\pi(t)]=\sigma}\ind_{H[\pi(1),\ldots,\pi(r)]\cong F}\ind_{H[\pi(1),\ldots,\pi(t),\pi(r+1),\ldots,\pi(6)]\cong G}.
 \label{eq:ordinary-M}
\end{align}
In each induced subgraph in \eqref{eq:ordinary-M}, the vertex $\pi(i)$ carries label $i$ for $1\le i\le t$. The other vertices are unlabeled, and flag isomorphisms fix each label. For positive semidefinite $Q_\sigma$, the evaluation of $\sum_H\tr(Q_\sigma M_\sigma(H))H$ is nonnegative up to $O(n^{-1})$. Conditioning on the type reduces it to an average of quadratic forms of independent extension densities.

For an unlabeled graph $K$ on four vertices, let
\begin{align}
 D_K(H)&=\frac{1}{6!}\sum_{\pi\in S_6}
       \ind_{H[\pi(1),\ldots,\pi(4)]\cong K},\notag\\
 E_K(H)&=\frac{1}{6!}\sum_{\pi\in S_6}
       \ind_{H[\pi(1),\ldots,\pi(4)]\cong K}
                       a_{\pi(5)\pi(6)}.
 \label{eq:ordinary-products}
\end{align}
Their evaluations are $p(K,Q)+O(n^{-1})$ and $sp(K,Q)+O(n^{-1})$, where $s=\frac{2e(Q)}{n^2}$. On an interval $s\in[s_-,s_+]$, the two expressions $E_K-s_-D_K$ and $s_+D_K-E_K$ have nonnegative limiting values. There are $11$ possibilities for $K$.

For rational $\iota,\rho,\beta,s_-,s_+$, the finite certificate condition is
\begin{align}
 f(H)-\iota-\rho\frac{e(H)}{15}
 &-\sum_K\lambda_K^-(E_K(H)-s_-D_K(H))\notag\\
 &-\sum_K\lambda_K^+(s_+D_K(H)-E_K(H))\notag\\
 &-\sum_\sigma\tr(Q_\sigma M_\sigma(H))\ge\beta
                      \qquad(H\in\mathcal G_6),
 \label{eq:ordinary-certificate}
\end{align}
with $\lambda_K^\pm\ge0$ and $Q_\sigma\succeq0$. Since the six-vertex densities sum to $1$, this condition gives
\begin{align}
 \mathcal F(Q)\ge\iota+\rho s+\beta-o(1)
          \qquad(s\in[s_-,s_+]).\label{eq:ordinary-bound}
\end{align}

The same estimate holds when the density is within $o(1)$ of the stated interval and the affine bound is evaluated at the nearest point of the interval. Indeed, the interval expressions then have evaluations bounded below by $-o(1)$, and their multipliers are fixed. This also accounts for replacing $\frac{2e(Q)}{M^2}$ by $1-c$ in the padded model.

\subsection{The additional objectives}

For the rooted terms in \eqref{eq:moment-affine}, let $\bar a_{ij}=1-a_{ij}$ for distinct vertices $i,j$ of a graph $H$ of order six. Let $x,y,z,p,q,r$ be a permutation of the vertices of $H$. Define
\begin{align*}
 T_{xyz}&=\bar a_{xy}\bar a_{xz}\bar a_{yz},\\
 S_{xyzp}&=\bar a_{xy}\bar a_{xz}\bar a_{xp},\\
 K_{xyzp}&=S_{xyzp}\bar a_{yz}\bar a_{yp}\bar a_{zp}.
\end{align*}
Averaging $T_{xyz}$ and $T_{xyz}T_{xpq}$ over all permutations gives the coefficients of $\langle a\rangle$ and $\langle a^2\rangle$, respectively. The averages of $S_{xyzp}$ and $S_{xyzp}T_{xqr}$ give the coefficients of $\langle r^3\rangle$ and $\langle r^3a\rangle$. The averages of $K_{xyzp}$ and $K_{xyzp}T_{xqr}$ give the coefficients of $\langle b\rangle$ and $\langle ba\rangle$. In each case, the common root is $x$, and unused vertices of the permutation can be ignored. Replacing independent samples by distinct samples changes the averages by $O(n^{-1})$. Thus \eqref{eq:ordinary-certificate} also applies with $f(H)$ increased by the appropriate linear combination of these coefficients.

More generally, for a monomial $r^ia^jb^k$ appearing after expansion of \eqref{eq:mixed-moment-affine}, use one common root and $i$ groups of one further vertex, $j$ groups of two further vertices, and $k$ groups of three further vertices. For each group, require all edges in the complement on the root and the vertices of that group. The coefficient is the average of the product of the corresponding $\bar a_{uv}$ over all bijections of the six vertices. Averages of empty products are $1$. Since $1+i+2j+3k\le6$ for all these monomials, the same coefficient test applies.

\subsection{The certified inequalities}

In Table~\ref{tab:ordinary}, each bonus $\beta_i$ gives the bound $\mathcal F(Q)\ge0.018+0.385s+\beta_i-o(1)$ on its interval. Let $s_1,\ldots,s_{18}$ be the left endpoints in increasing order, and let $s_{19}=+\infty$ and $\eta=10^{-4}$. For $s\in[s_1,s_2)$, let $\beta(s)=\beta_1$. For $2\le i\le18$ and $s\in[s_i,s_{i+1})$, let
\begin{align}
 \beta(s)=\beta_{i-1}+(\beta_i-\beta_{i-1})
                    \min\left\{1,\frac{s-s_i}{\eta}\right\}.
 \label{eq:bonus}
\end{align}
The bonuses are nondecreasing, so this interpolation decreases the bound on each interval where a bonus changes.

\begin{table}[htbp]
\centering
\caption{Bonuses for the affine bound on $\mathcal F(Q)$.}
\label{tab:ordinary}
\begin{tabular}{@{}rrr@{\qquad}rrr@{}}
\toprule
$s_-$&$s_+$&$\beta_i$&$s_-$&$s_+$&$\beta_i$\\
\midrule
0.250&0.255&0.00284 & 0.310&0.320&0.00545\\
0.255&0.260&0.00311 & 0.320&0.330&0.00579\\
0.260&0.265&0.00334 & 0.330&0.340&0.00611\\
0.265&0.270&0.00358 & 0.340&0.355&0.00637\\
0.270&0.275&0.00382 & 0.355&0.370&0.00675\\
0.275&0.280&0.00404 & 0.370&0.390&0.00701\\
0.280&0.290&0.00426 & 0.390&0.410&0.00732\\
0.290&0.300&0.00468 & 0.410&0.430&0.00741\\
0.300&0.310&0.00508 & 0.430&0.460&0.00741\\
\bottomrule
\end{tabular}
\end{table}

Table~\ref{tab:tail-gains} gives two further families of bounds. For $j=1,2$, let $T_j(s)$ be the indicated function on $[0.43,0.89]$, with $T_j(s)=0$ for $s<0.43$. Extend each function by its last constant value for $s>0.89$. The affine formulas agree at their shared endpoints, and both functions are nondecreasing.

\begin{table}[htbp]
\centering
\caption{Further affine bounds on $\mathcal F(Q)$.}
\label{tab:tail-gains}
\begin{tabular}{@{}cll@{}}
\toprule
Function & Interval & Lower bound\\
\midrule
$T_1$ & $[0.43,0.53]$ & $0.02835+0.375s$\\
      & $[0.53,0.62]$ & $0.05379+0.327s$\\
      & $[0.62,0.70]$ & $0.06681+0.306s$\\
      & $[0.70,0.77]$ & $0.08851+0.275s$\\
      & $[0.77,0.89]$ & $0.30026$\\
\midrule
$T_2$ & $[0.43,0.50]$ & $0.02750+0.380s$\\
      & $[0.50,0.54]$ & $0.03700+0.361s$\\
      & $[0.54,0.60]$ & $0.05644+0.325s$\\
      & $[0.60,0.66]$ & $0.06364+0.313s$\\
      & $[0.66,0.70]$ & $0.07156+0.301s$\\
      & $[0.70,0.77]$ & $0.08976+0.275s$\\
      & $[0.77,0.89]$ & $0.30151$\\
\bottomrule
\end{tabular}
\end{table}

\begin{table}[htbp]
\centering
\caption{Bounds on $\mathcal F(Q)+\langle\omega_v(a_v-U_*)\rangle$.}
\label{tab:moment-gains}
\begin{tabular}{@{}lcll@{}}
\toprule
Weight $\omega_v$ & $U_*$ & Interval & $L_*(s)$\\
\midrule
$0.1a_v$ & $0.617$ & $[0.255,0.48]$ & $0.03775+0.351s$\\
\addlinespace
$0.18r_v^3+0.08b_v$ & $0.618$ & $[0.25,0.30]$ & $0.06410+0.281s$\\
                   &          & $[0.30,0.42]$ & $0.06140+0.290s$\\
\addlinespace
$0.24r_v^3+0.025b_v$ & $0.618$ & $[0.25,0.32]$ & $0.06865+0.268s$\\
                    &          & $[0.32,0.42]$ & $0.06321+0.285s$\\
\addlinespace
$0.24r_v^3+0.025b_v$ & $0.558$ & $[0.25,0.32]$ & $0.07600+0.268s$\\
                    &          & $[0.32,0.42]$ & $0.07760+0.263s$\\
\bottomrule
\end{tabular}
\end{table}

\begin{proposition}\label{prop:ordinary-finite}
The bounds in Tables~\ref{tab:ordinary}--\ref{tab:moment-gains} and the $33$ bounds of the form \eqref{eq:mixed-moment-affine} supplied in \cite{certificates} hold with uniform error $o(1)$. We also have $\mathcal F(Q)\ge0.20303-o(1)$ on $[0.46,0.525]$. Consequently the probability of no offer among the $M$ candidates is at least $L_0(s)-o(1)$ for $s\in[0.25,1]$.
\end{proposition}

\begin{proof}
The rational witnesses in \cite{certificates} give \eqref{eq:ordinary-certificate} on each interval, or on a subdivision of that interval with the same affine bound. For the additional objectives, use the coefficients in the preceding subsection. For each witness, the $22$ interval multipliers are nonnegative, each of the $14$ matrices has positive pivots under the elimination recurrence \eqref{eq:bareiss}, and the coefficient inequality holds on every $H\in\mathcal G_6$. The verification constructs these coefficients by summing over the vertex permutations and checks the inequalities with exact arithmetic. For each subdivision, the interval endpoints are also checked to cover the stated interval. In \texttt{moments.json}, each row gives the endpoints $s_-,s_+$, the parameters $U_*,R_*$, the slope $\rho$, and the intercept $\iota$ in units of $10^{-8}$. The two arrays of weights give the coefficients of $\omega_v$ and $\zeta_v$ in the order in \eqref{eq:mixed-monomials}, in the same units. For these rows, take $\beta=0$ in \eqref{eq:ordinary-certificate}.

These inequalities and \eqref{eq:no-offer-F} give the first expression inside the maximum in \eqref{eq:L0} on $[0.25,0.525]$. Indeed, truncating the linear part at $s=0.446$ and interpolating the bonuses only decrease the affine bound. On $[0.46,0.525]$, the bound $0.20303$ exceeds $0.19712$, the value of the expression there. Applying the proof of the second inequality in \eqref{eq:three-term} to all $M$ padded candidates with $t=1$ gives $\Prb(Z=0)\ge g(s)-o(1)$ throughout $[0.25,1]$. Since $g$ is increasing and $g(0.525)>0.19712$, this also gives the first term in \eqref{eq:L0} for $s\ge0.525$. On $[0.43,0.89]$, Table~\ref{tab:tail-gains} and \eqref{eq:no-offer-F} give the bounds $T_1(s)$ and $T_2(s)$, up to $o(1)$. Their constant extensions for $s>0.89$ are bounded above by $g(s)$, since $g(0.89)>0.30151$. For $s<0.43$, both functions are zero. Taking the maximum gives $\Prb(Z=0)\ge L_0(s)-o(1)$.
\end{proof}

\section{The palette calculation}\label{app:palette}

Let $M_0=10^8$ and $S_0=10^{12}$ be the arithmetic scales for weights and probabilities, and for palette costs and scales, respectively. An integer $W$ representing a weight or probability has value $\frac{W}{M_0}$, and an integer $Z$ representing a palette cost or scale has value $\frac{Z}{S_0}$. Let $L_i=M_0\ell_i$ for $1\le i\le176$, as given in Table~\ref{tab:cutoffs}, and let $L_{177}=M_0$. We give the rectangle bounds and the recurrence used in Proposition~\ref{prop:palette}.

\begin{table}[htbp]
\centering
\caption{Degree cutoffs $L_1,\ldots,L_{176}$, read by rows.}\label{tab:cutoffs}
\small
\setlength{\tabcolsep}{4pt}
\begin{tabular}{@{}rrrrrrrr@{}}
\toprule
81750000 & 81751280 & 81754612 & 81759766 & 81766629 & 81775128 & 81785208 & 81796827\\
81809949 & 81824544 & 81840586 & 81858054 & 81876926 & 81897184 & 81918811 & 81941793\\
81966116 & 81991766 & 82018732 & 82047002 & 82076565 & 82107412 & 82139535 & 82172923\\
82207567 & 82243462 & 82280598 & 82318968 & 82358565 & 82399383 & 82441415 & 82484655\\
82529098 & 82574737 & 82621566 & 82669582 & 82718777 & 82769149 & 82820691 & 82873399\\
82927268 & 82982295 & 83038473 & 83095802 & 83154274 & 83213888 & 83274637 & 83336519\\
83399532 & 83463670 & 83528930 & 83595309 & 83662805 & 83731412 & 83801129 & 83871951\\
83943877 & 84016903 & 84091026 & 84166243 & 84242552 & 84319950 & 84398434 & 84478001\\
84558648 & 84640375 & 84723177 & 84807052 & 84891999 & 84978013 & 85065094 & 85153239\\
85242444 & 85332711 & 85424034 & 85516412 & 85609842 & 85704324 & 85799854 & 85896432\\
85994054 & 86092719 & 86192424 & 86293169 & 86394950 & 86497767 & 86601618 & 86706499\\
86812410 & 86919350 & 87027316 & 87136306 & 87246319 & 87357353 & 87469406 & 87582478\\
87696566 & 87811668 & 87927783 & 88044910 & 88163047 & 88282193 & 88402345 & 88523503\\
88645665 & 88768829 & 88892994 & 89018159 & 89144322 & 89271483 & 89399638 & 89528788\\
89658931 & 89790065 & 89922189 & 90055302 & 90189403 & 90324489 & 90460562 & 90597617\\
90735655 & 90874674 & 91014674 & 91155651 & 91297608 & 91440540 & 91584447 & 91729328\\
91875183 & 92022009 & 92169806 & 92318573 & 92468308 & 92619010 & 92770677 & 92923312\\
93076909 & 93231470 & 93386992 & 93543475 & 93700919 & 93859321 & 94018681 & 94178997\\
94340270 & 94502497 & 94665678 & 94829812 & 94994898 & 95160934 & 95327920 & 95495856\\
95664739 & 95834570 & 96005346 & 96177067 & 96349734 & 96523343 & 96697895 & 96873389\\
97049823 & 97227197 & 97405510 & 97584760 & 97764949 & 97946074 & 98128133 & 98311128\\
98495056 & 98679917 & 98865711 & 99052436 & 99240090 & 99428675 & 99618190 & 99808631\\
\bottomrule
\end{tabular}
\end{table}

\subsection{The slopes in the density estimate}

Let $W_1\le\cdots\le W_m=M_0$ be the current upper weights. Take the lower convex hull of the points $(W_i,L_i)$. If several points have the same first coordinate, retain only the one with the smallest second coordinate. Remove collinear interior hull vertices. For each consecutive pair $(W_i,L_i),(W_j,L_j)$ on the hull, include the pair
\begin{align}
 (p,q)=(L_i-L_j,W_j-W_i).\label{eq:slope-pair}
\end{align}
Its slope parameter is $\lambda=\frac{p}{q}\le0$. Also include $(0,1)$. If there is an index with $L_i<86485000$ and $W_i<M_0$, include a pair representing
\begin{align}
 \min_{\substack{L_i<86485000\\W_i<M_0}}
       \frac{L_i-86485000}{M_0-W_i}.
 \label{eq:last-slope}
\end{align}
Let $\mathcal A$ be the resulting finite collection of pairs. For $(p,q)\in\mathcal A$, let
\begin{align}
 A=\max_i\{q(86485000-L_i)-pW_i\},\quad
 B=\max\{0,-pM_0-A\}.\label{eq:integer-AB}
\end{align}
These are integer versions of $A_\lambda,B_\lambda$ in \eqref{eq:AB}. Both $\lambda$ and the value $J(d)$ defined in \eqref{eq:rectangle-J} depend on $p,q,A,B$ only through $\frac{p}{q}$, $\frac{A}{q}$, and $\frac{B}{q}$.

\subsection{An upper bound on a rectangle}

Let
\begin{align}
 \mathcal D=[d_-,d_+]\times[u_-,u_+]
       \subseteq[L_i,L_{i+1}]\times[0,M_0]
 \label{eq:integer-rectangle}
\end{align}
be a rectangle with integer endpoints. Its coordinates represent $M_0d$ and $M_0u$ in \eqref{eq:u}. Let $W=W_1$ and let
\begin{align}
 \mu_-&=W+\floor{\frac{(M_0-W)u_-}{M_0}},\quad
 \mu_+=W+\ceil{\frac{(M_0-W)u_+}{M_0}},\notag\\
 t_+&=\ceil{\frac{d_+\mu_+}{M_0}},\quad
 V=\ceil{\frac{d_+^2}{M_0}}.
 \label{eq:rectangle-t}
\end{align}
The integer $t_+$ bounds $M_0T$ from above.

For a pair $(p,q)$ and its $A,B$ in \eqref{eq:integer-AB}, let
\begin{align}
 N(d)&=M_0B+d(p\mu_-+A-B),\notag\\
 K(d)&=
 \begin{cases}
 185481000,&N(d)\ge0,\\
 185480000,&N(d)<0,
 \end{cases}\notag\\
 J(d)&=\ceil{\frac{K(d)N(d)}{qM_0^2}}.
 \label{eq:rectangle-J}
\end{align}
The two choices for $K(d)$ use the bounds $1.8548<k<1.85481$ with the correct direction for the sign of $N(d)$. In particular, $J(d)$ is an upper bound for the signed degree term multiplied by $kM_0$.

For each $(\gamma,b_\gamma)\in\mathcal P$, let
\begin{align*}
 D_\gamma(d)=M_0(b_\gamma-\gamma)+\ceil{\gamma d}.
\end{align*}
For each $(p,q)\in\mathcal A$ and $(\gamma,b_\gamma)\in\mathcal P$, replace $V$ by the minimum of its current value and
\begin{align*}
 \max\{D_\gamma(d_-)+J(d_-),\,D_\gamma(d_+)+J(d_+)\}.
\end{align*}
For each $(\kappa,\gamma,b')\in\mathcal P'$, do the same with $D_\gamma(d)$ replaced by $M_0(b'-\gamma)+\ceil{\gamma d}$. In $J(d)$, use $\ceil{M_0\kappa}$ for nonnegative $N(d)$ and $\floor{M_0\kappa}$ for negative $N(d)$. If $\kappa=0$, only the affine part remains, so a single endpoint comparison suffices. For every real pair $(d,u)$ with $(M_0d,M_0u)\in\mathcal D$, we have $[C(d,T)]_+\le\frac{V}{M_0}$. Indeed, since $\lambda\le0$, replacing the mean weight $w+(1-w)u$ by $\frac{\mu_-}{M_0}$ gives an upper bound on $H_\lambda(d,T)$. At that mean weight each expression in \eqref{eq:C}, other than $d^2$, is affine in the degree, so its maximum is at one of the two endpoints. Thus $C(d,T)\le\frac{V}{M_0}$, and the verifier checks $V\ge0$.

Let $\bar u$ be the point of $[u_-,u_+]$ nearest to $\frac{M_0}{2}$, and let $H_0=\ceil{\frac{d_+u_+}{M_0}}$ and $W^{(2)}=\ceil{\frac{W^2}{M_0}}$. The upper capacities, the allocation, and the resulting bounds on the density are
\begin{align}
 A_+&=\ceil{\frac{H_0^2}{M_0}},\quad
 B_+=\ceil{\frac{2d_+^2\bar u(M_0-\bar u)}{M_0^3}},\notag\\
 h&=\min\{V,A_+\},\quad
 j=\min\{V-h,B_+\},\notag\\
 C_+&=\ceil{\frac{M_0 h+Wj+W^{(2)}(V-h-j)}{M_0}},\quad
 c_+=\min\left\{C_+,\ceil{\frac{t_+^2}{M_0}}\right\}.
 \label{eq:integer-density}
\end{align}
Here and below a floor or ceiling is taken after exact rational evaluation of its argument. In \eqref{eq:integer-density}, $C_+$ and $c_+$ represent the corresponding bounds from Section~\ref{sec:profiles} multiplied by $M_0$. The verifier checks $M_0-c_+\ge25000000$.

For $G_+$ in \eqref{eq:mu-budget}, evaluate
\begin{align*}
 D_\gamma(d)+\ceil{\frac{185481000\bigl((86485000-L_1)d+(M_0-86485000)(M_0-d)\bigr)}{M_0^2}}
\end{align*}
at $d_-,d_+$ for every pair in $\mathcal P$. Take the maximum over the two endpoints and then the minimum over the pairs. The resulting integer is an upper bound on $M_0G_*(d)$ throughout the rectangle. For each triple in $\mathcal P'$ with $\kappa>0$, use the analogous endpoint maximum for $M_0G_+^{\kappa,\gamma,b'}$, with the corresponding $b',\gamma$ and with $\ceil{M_0\kappa}$ in place of $185481000$.

\subsection{Bounds on the rooted triangle counts}

In this subsection and the next, expressions in real units are evaluated on the grid $\frac{1}{M_0}\mathbb Z$. Rounding up and down mean taking the least grid point at or above the value and the greatest grid point at or below the value, respectively.

For \eqref{eq:initial-triangle-cap}, use $1.85481$ in place of $k$, round $k(a-\tau)$ up first, and round the resulting $U_\gamma$ up. Use the pair with the smallest resulting $U_\gamma$ from Table~\ref{tab:root}, using the first such pair in case of a tie. For \eqref{eq:retention-triangle-cap}, round $\beta_w$ and then $(1+w)\beta_w$ up, and check the condition $(1+w)\beta_w\le1.8548w^2$ with these upper bounds. If the condition holds, round both terms of $U$ and the coefficient $\nu$ up and include the resulting pair.

On an integer rectangle $\mathcal D$, the quantity $c_+-\floor{\nu(d_--L_1)}$ bounds $M_0\bigl(c^*(d,u)-\nu(d-\tau)\bigr)$ for every real pair $(d,u)$ with $(M_0d,M_0u)\in\mathcal D$. For each of the nine slopes $\nu$ in Section~\ref{sec:profiles}, the search applies the subdivision rule below to $[L_1,M_0]\times[0,M_0]$ with tolerance $4200$. The initial point rectangles have degree in $\{L_1,M_0\}$ and second coordinate in $\{0,\frac{M_0}{4},\frac{M_0}{2}, \frac{3M_0}{4},M_0\}$. The resulting upper bounds, divided by $M_0$, give nine more values of $U$ with their respective slopes $\nu$.

\subsection{The offer bound on a rectangle}

Use the real units $\frac{t_+}{M_0}$, $\frac{C_+}{M_0}$, and $\frac{c_+}{M_0}$ for the bounds in Section~\ref{sec:triangle-offer}. Evaluate the functions in \eqref{eq:L0} at $1-\frac{c_+}{M_0}$. Round the affine terms and the bonus interpolation down. For any square root $\sqrt{z}$ needed in these bounds, the verifier uses the upper bound $\frac{\floor{M_0\sqrt{z}}+1}{M_0}$. Thus it uses the grid point immediately above the floor, including when $\sqrt{z}$ is on the grid. The floor is verified by integer comparisons of squares. In $g$, use this upper square root and round the resulting value of $g$ down.

For each of the pairs $(U,\nu)$ just obtained, apply the interval and slope tests of Section~\ref{sec:triangle-offer} at $s=1-\frac{C_+}{M_0}$. Use $1.8548$ as the lower bound on $k$ in \eqref{eq:mu-budget} and \eqref{eq:slope-test}. For the weight $0.1a_v$, use the exact value $\phi=\frac{0.1C_+}{M_0}$. For the star weights, round both summands of $\phi$ in \eqref{eq:star-weight-bounds} up, using the upper bound on the square root. In each case, round the derivative term $g_0A_0'\bigl(\frac{C_+}{M_0}\bigr)$ up in the slope test. Evaluate the test after clearing denominators. For \eqref{eq:moment-lower}, evaluate the applicable affine bounds at $s=1-\frac{C_+}{M_0}$ and round each value down to the grid. For the weight $0.24r_v^3+0.025b_v$, interpolate between the rounded values and round the result down again when interpolation is used. Round the terms $\nu\phi m_0$ and $W_0h_0$ up, and round the positive correction $g_0A_0\!\left(\frac{C_+}{M_0}\right)$ down. For the star weights, round $\left(\frac{C_+}{M_0}\right)^3$ down to the $\frac{1}{M_0}$ grid before multiplying by $g_0\eta_3$.

For \eqref{eq:mixed-moment-lower}, use the same reference density $\frac{C_+}{M_0}$, before taking the minimum with $t_+^2$ in real units. Round pointwise and average upper bounds up and average lower bounds down, including the intermediate products. To compute a lower bound on $A$ in \eqref{eq:triangle-average-support} at the reference density, round $\frac{c^2}{t_+}$ down and $t_+c$ up. The resulting lower bound differs from $A(c)$ by at most $\frac{3}{M_0}$. Add $\frac{3}{M_0}$ to it when bounding the derivatives of the lower bounds for $\langle ra\rangle$, $\langle r^2a\rangle$, and $\langle a^2\rangle$ from above. For all derivative terms, use upward rounding, including at the positive part in \eqref{eq:triangle-average-support}. In the bound using $\xi$, round the first term in each difference up and the second term down.

For the affine upper bounds on $\Gamma\mu$, round both $p$ and $z$ up and use the lower grid bounds on the positive signed coefficients in the denominators. Let $z_0$ be the slope $\rho$ minus the upper bound on $A_*'(C_+)$ in real units. If $z_0\ge0$, take the smallest intercept among the bounds with $z\le z_0$ and the convex combinations of two bounds with slopes on opposite sides of $z_0$ evaluated at slope $z_0$. Round each interpolated intercept up. In \eqref{eq:mixed-moment-lower}, round all favorable contributions down and all subtracted contributions up. Take the maximum of the capped lower bounds from \eqref{eq:moment-lower} and \eqref{eq:mixed-moment-lower} and the lower bound from \eqref{eq:L0}.

Let $x=1-\frac{t_+}{M_0}$. In the padded bound in \eqref{eq:rectangle-offer}, use
\begin{align*}
 e^x\ge 1+x+\frac{x^2}{2}+\frac{x^3}{6}+\frac{x^4}{24}.
\end{align*}
Round each term down, and round down the product with the no-offer bound before subtracting it from $1$. For the direct bound in \eqref{eq:rectangle-offer}, round $t_+^2$ down in real units before subtracting the upper density bound and taking the positive part. Round the square root according to the same convention and the value of $g$ down. The smaller of the two resulting bounds, multiplied by $M_0$, is the integer upper bound used for the rectangle.

\subsection{Subdivision of rectangles}

For class $i$, the initial rectangle is $[L_i,L_{i+1}]\times[0,M_0]$. If $W=M_0$, the upper bounds depend only on the first coordinate, and we may use the interval $[0,0]$ for the second coordinate. Let $e$ be the largest upper bound on the point rectangles with degree in $\{L_i,L_{i+1}\}$ and second coordinate in $\{0,\frac{M_0}{10},\frac{M_0}{5},\frac{M_0}{2}, \frac{9M_0}{10},M_0\}$. Let $\varepsilon_0$ be the integer search tolerance for the current phase.

For each rectangle, let $p$ be its upper bound. If $p\le e+\varepsilon_0$, finish processing the rectangle. Otherwise evaluate the bound at its midpoint, with each coordinate rounded down to an integer, and increase $e$ to that value if it is larger. Test $p\le e+\varepsilon_0$ again. If the inequality fails and both widths are at most $1$, increase $e$ to $p$. In the remaining case, bisect at the midpoint of the degree interval if its width is greater than $1$ and either the second width is at most $1$ or
\begin{align*}
 d_+-d_-\ge \floor{\frac{(M_0-W)(u_+-u_-)}{M_0}}.
\end{align*}
Otherwise bisect the second interval. Process the half with larger coordinates first. The procedure terminates since each split reduces a width, and the rectangles cover the initial domain. Every processed rectangle has upper bound at most the final $e+\varepsilon_0$ or has been covered by smaller rectangles. Thus $e+\varepsilon_0$ is an upper bound on the initial domain. For the bounds in \eqref{eq:scan-triangle-cap}, use the same rule with the corresponding initial domain, point rectangles, and tolerance. The program \texttt{palette.cpp} generates \texttt{palette\_bounds.txt} by this procedure, with tolerance $5800$ for the completion step. The verifier recomputes the analytic triangle bounds and checks the numerical triangle and offer-probability bounds by subdivision. It uses~\eqref{eq:integer-recurrence} and~\eqref{eq:integer-tail} to recompute the palette total.

\subsection{The recurrence and its value}

Let $W_i^{(j)}$ be the integer upper relative weight of class $i$ in phase $j$, let $Z_j$ be the integer upper scale, and let $A_j$ be the accumulated integer upper bound on $S_0$ times the sum of the normalized leading costs of the preceding phases. Initially $W_i^{(0)}=M_0$ for every $i$, $Z_0=S_0$, and $A_0=0$. Continue while $Z_j>2200000000$. The phase parameter $q_j$ and the search tolerance $\varepsilon_0$ depend on the rounded upper scale as follows.

\begin{table}[htbp]
\centering
\caption{Phase parameters.}\label{tab:phase-parameters}
\begin{tabular}{@{}lrr@{}}
\toprule
Upper scale $\frac{Z_j}{S_0}$ & $q_j$ & $\varepsilon_0$\\
\midrule
$[0.4,1]$ & 110 & 5800\\
$[0.1,0.4)$ & 110 & 5800\\
$[0.02,0.1)$ & 85 & 11600\\
$(0,0.02)$ & 45 & 23200\\
\bottomrule
\end{tabular}
\end{table}

Let $\pi_i$ be the maximum of the integer bounds from the subdivision for classes $1,\ldots,i$, and let $\pi_*=\pi_m$. The verifier checks $0<\pi_i\le\pi_*<2\pi_i$. Let $u_i=\frac{\pi_*-\pi_i}{\pi_i}$. By Taylor's formula on $[0,\frac{1}{q_j}]$, for $0\le u_i<1$ we have
\begin{align}
 \left(1-\frac{1}{q_j}\right)^{u_i}
 \le 1-\frac{u_i}{q_j}-\frac{u_i(1-u_i)}{2q_j^2}.
 \label{eq:weight-taylor}
\end{align}
Indeed, the third derivative of $(1-x)^{u_i}$ is nonpositive for $0\le x<1$. The updates are
\begin{align}
 A_{j+1}&=A_j+\ceil{\frac{Z_j\pi_*}{M_0q_j}},\notag\\
 W_i^{(j+1)}
   &=\ceil{W_i^{(j)}
      \left(1-\frac{u_i}{q_j}-\frac{u_i(1-u_i)}{2q_j^2}\right)},\notag\\
 Z_{j+1}&=\ceil{\frac{(q_j-1)Z_j}{q_j}}.
 \label{eq:integer-recurrence}
\end{align}
By \eqref{eq:trajectory} and \eqref{eq:weight-taylor}, the new upper weights dominate the actual relative retentions. The new upper scale dominates the actual scale. Since $Z_j\ge S_0R_j$ and $\pi_*\ge M_0p_{*,j}$, the first line adds an integer upper bound on $\frac{S_0p_{*,j}R_j}{q_j}$.

At termination, verify the recorded completion bounds using the final weights, and let $\pi_0$ be their maximum. Let $j$ be the number of phases. Since $Z_j\ge S_0R_{\mathrm{end}}$ and $\pi_0\ge M_0p_0$, the final update adds an integer upper bound on $S_0p_0R_{\mathrm{end}}$. Let
\begin{align}
 A_{\mathrm{final}}=A_j+\ceil{\frac{Z_j\pi_0}{M_0}}.
 \label{eq:integer-tail}
\end{align}
Running the verifier in \cite{certificates} gives
\begin{align}
 j=487,\qquad A_{\mathrm{final}}=819392866996.
 \label{eq:integer-output}
\end{align}

The file \texttt{palette\_bounds.txt} contains one row for each of the $487$ phases and one row for the completion step. At every phase, the upper weights are positive and nondecreasing in the class index, with the last weight equal to $M_0$. Together with the rounding rules, \eqref{eq:integer-output} proves \eqref{eq:palette-total}.

\end{document}